\documentclass{amsart}
\usepackage{graphicx,amssymb,epsfig}
\newtheorem{theorem}{Theorem}[section]
\newtheorem{proposition}[theorem]{Proposition}
\newtheorem{lemma}[theorem]{Lemma}
\newtheorem{corollary}[theorem]{Corollary}
\newtheorem{remark}[theorem]{Remark}
\newtheorem{example}[theorem]{Example}

\newtheorem{definition}[theorem]{Definition}

\newcommand{\bth}{\begin{theorem}}
\newcommand{\bpr}{\begin{proposition}}
\newcommand{\epr}{\end{proposition}}
\newcommand{\bco}{\begin{corollary}}
\newcommand{\eco}{\end{corollary}}
\newcommand{\ble}{\begin{lemma}}
\newcommand{\ele}{\end{lemma}}
\newcommand{\bre}{\begin{remark}\rm}
\newcommand{\ere}{\end{remark}}
\newcommand{\bex}{\begin{example}\rm}
\newcommand{\eex}{\end{example}}
\newcommand{\bde}{\begin{definition}\rm}
\newcommand{\ede}{\end{definition}}

\usepackage[all]{xy}

\def\la#1{\hbox to #1pc{\leftarrowfill}}
\def\ra#1{\hbox to #1pc{\rightarrowfill}}
\def\fract#1#2{\raise3pt\hbox{$ #1 \atop #2 $}}

\def\lrar{{\ra 2}}

\def\tensor{\otimes}
\def\sp#1{\hbox{SP}^{#1}}
\def\sj#1{\hbox{Sym}^{* {#1}}}
\def\bsp#1{\overline{\hbox{SP}}^{#1}}

\def\spy{\hbox{SP}^{\infty}}

\def\bbz{{\mathbb Z}}
\def\bbf{{\mathbb F}}
\def\bbp{{\mathbb P}}
\def\bbr{{\mathbb R}}

\def\bbr{{\mathbb R}}

\font\sc=cmcsc10 at 11pt

\def\lrar{{\ra 2}}

\begin{document}

\title{Symmetric Joins and Weighted Barycenters}
\author{Sadok Kallel, Rym Karoui}
\address{Universit\'e des Sciences et Technologies de Lille\\
 Laboratoire Painlev\'e, U.F.R de Math\'ematiques\\
59655 Villeneuve d'Ascq, France}
\email{sadok.kallel@math.univ-lille1.fr}

\maketitle

\begin{abstract}
Given a space $X$, we study
    the homotopy type of ${\mathcal B}_n(X)$ the space obtained as the ``union
    of all $(n-1)$-simplexes spanned by points in $X$" or the space of
    ``formal barycenters of weight $n$ or less" of $X$. This is a space
    encountered in non-linear analysis under the name of \textit{space of
      barycenters} or in differential geometry in the case $n=2$ as the
    \textit{space of chords}. We first relate this space to a more familiar
    symmetric join construction and then determine its stable homotopy type in
    terms of the symmetric products on the suspension of $X$. This leads to a
    complete understanding of the homology of ${\mathcal B}_n(X )$ as a
    functor of $X$, and to an expression for its Euler characteristic given in
    terms of that of $X$.  A sharp connectivity theorem is also established.
    Finally the case of spheres $S$ is studied in details and the homotopy
    type of $\mathcal B_n(S)$ is described generalizing in this way an early
    and beautiful result of James, Thomas, Toda and Whitehead.
\end{abstract}

%***********************************************************************

\section{Introduction}\label{intro}

An interesting construction of Bahri and Coron \cite{bahri} associates to a
given topological space $X$ and integer $n > 0$, the space ${\mathcal B}_n(X)$
of ``weighted barycenters" of $X$ obtained by taking the union of all
$(n-1)$-simplexes $\Delta (p_1,\ldots, p_n)$ on vertices $p_i\in X$ with the
topology that when two vertices of $\Delta$ come together, the resulting
$n-2$-simplex is identified with a face of $\Delta$ (precise definitions in
\S\ref{barycenters}).  In the notation of \cite{bahri} these spaces were
described according to
$${\mathcal B}_n(X) = \left\{\sum_{i=1}^n t_i\delta_{x_i}\ |\ (x_1,\ldots, x_n)\in X^n,
(t_1,\ldots ,t_n)\in\Delta_{n-1}\right\}$$
where $\delta_x$ is the ``Dirac
mass" at the point $x$ and $\Delta_{n-1}$ is the $(n-1)$-simplex of
all tuples $(t_1,\ldots, t_n)$, $t_i\geq 0$ and $\sum t_i=1$. The formal
sums are understood to be \textit{abelian}; i.e. unordered, and if two $x$'s
coincide their coefficients add up. An alternative geometric description
can also be given as follows. Let $X$ be a finite CW
complex. Then there is some choice of embedding $i:
X\hookrightarrow\bbr^N$, $N$ large, whereby the space ${\mathcal B}_n(X)$ is the
union of all $(n-1)$-dimensional simplices with vertices in $i(X)$
such that the intersection of any two such simplices is a common
face if any (see \cite{gorinov}).

Spaces of barycenters play an important role in the study of limiting Sobolev
exponent problems in non-linear analysis such as the Yamabe and the
scalar-curvature equations. Symmetric products on the other hand (see below)
appear in the study of \textit{singular} solutions to these equations. The
smooth solutions come as critical points of a functional, with nested level
sets $\cdots W_p\subset W_{p+1}\subset\cdots$, and when no critical point
exists, the pair $(W_{p+1},W_p)$ behaves topologically much like $({\mathcal
  B}_{p+1}(M),{\mathcal B}_p(M))$, where $M$ is the underlying domain of the
equation. Induction formulae relating the $\bbz_2$-orientation class of the
pair $({\mathcal B}_{p+1}(M),{\mathcal B}_p(M))$ to that of $({\mathcal
  B}_{p}(M),{\mathcal B}_{p-1}(M))$ one level lower, are then used to derive
existence results for the solutions of such nonlinear elliptic equations (see
\cite{bahri}, with an appendix by Jean Lannes, or \cite{bb}). More recently,
\cite{dm} proved that barycenter spaces played a fundamental role in the
two-dimensional scalar curvature problem and other conformal equations.

The spaces ${\mathcal B}_n(X)$ enter differently but fundamentally as well in
work of Vassiliev as a tool to construct simplicial resolutions for
complements of discriminant loci in algebraic geometry \cite{gorinov,
  vassiliev}. This technique has been skillfully applied to provide stable
splittings for so-called ``Atiyah-Hitchin schemes" in \cite{casanave}.

The space ${\mathcal B}_2(X)$ is an even older concept. This is the so-called
space of \textit{chords} on $X$ and has been considered for instance in
\cite{clark} and much more recently in \cite{mineyev}. The terminology comes
from the fact that ${\mathcal B}_2(X)$ can be viewed as the space of all
``chords" $pq$, $p$ and $q$ in $X$, with the understanding that $pp$ is
identified with the point $p$ (see \S\ref{chords}).

The purpose of this paper is to investigate the homotopy type and homology of
${\mathcal B}_n(X)$ for all $n\geq 1$ and for $X$ any reasonable topological
space (i.e. based, connected, locally compact and hausdorff).  The nature of
these spaces is intimately related to that of more familiar constructions
known as symmetric joins and symmetric products.

We will write $X*Y$ the join product of $X$ and $Y$. The $n$-fold
\textit{symmetric join} of $X$ is obtained as the quotient of the
$n$-th iterated join $X^{*  n}:= X*  \cdots
*  X$ by the permutation action of the symmetric group $\mathfrak{S}_n$ and
is written $\sj{n}(X)$ (see \S\ref{constructions}).  This space
for $X$ a sphere and $n=2$ for example is treated in \cite{james}
where the following cute result was obtained (see corollary \ref{maincor})
\begin{equation}\label{jmea}
\sj{2}(S^k ) \cong \Sigma^{k+1}\bbr P^k
\end{equation}
Here $\Sigma^{k+1}X:=S^{k}*X$ refers to the $k$-fold {\sl unreduced}
suspension of $X$ and $\bbr P^k$ is the $k$-th real projective space.  Very
little else seems to be known about the homotopy type and homology of
symmetric joins and the purpose of this paper is to remedy to this.

The following is our first main observation.

\bth\label{main1} ${\mathcal B}_n(X)$ and $\sj{n}(X)$ have the same
homotopy type.
\end{theorem}

In particular when $X$ is a simplicial complex\footnote{in other terminology,
a \textit{polyhedral space}.},
we show that both spaces have the homotopy type of the quotient
$S^{n-1}\wedge_{\mathfrak{S}_n}X^{(n)}$,
where $X^{(n)}$ is the ``smash product" of $X$ (see \S\ref{constructions})
and where $\mathfrak{S}_n$ acts on $S^{n-1}\subset\bbr^n$ by
permutation of coordinates (details in \S\ref{barycenters}).
It is easy to see that $\sj{n}(X)$ has the structure
of a CW complex of dimension $n(d+1)-1$ if $X$ is a simplicial
complex of dimension $d$. This shows for instance that the homology of
${\mathcal B}_n(X)$ vanishes in degrees greater than $n(d+1)$. On the other
hand and when $X=M$ is a closed orientable manifold of dimension $d$, $\mathcal B_n(M)$
always has a top homology class mod-$2$ while with integer coefficients (see \S\ref{manifolds})
$$H_{n(d+1)-1}({\mathcal B}_n(M)) = \begin{cases}\bbz,& \hbox{if $d$ is odd}\\
0,& \hbox{if $d$ even}\end{cases}$$

Using a handy description of a symmetric join as an iterated pushout,
we show based on a Van-Kampen type argument that
$\sj{n}(X)$, and hence ${\mathcal B}_n(X)$, is simply connected as soon as $n> 1$
(Theorem \ref{simplyconnected}). The connectivity of barycenter spaces is however
more tricky to establish and we appeal to that end to some old results of Nakaoka
\cite{nakaoka} and a relative Leray-Serre spectral sequence argument (see \S6).
The optimal result is as follows.

\bth\label{cor1} If $X$ is an $r$-connected CW complex, $r\geq 1$,
then ${\mathcal B}_n(X)$ is $(2n+r-2)$-connected.
\end{theorem}

The above lower bound for the
connectivity is sharp in light of (\ref{jmea}).

It turns out that the best way to getting to the homology of $\sj{n}(X)$ or
${\mathcal B}_n(X)$ is to analyze their \textit{stable} homotopy type (i.e.
after suspension). Write $\sp{n}X$ for the $n$-th symmetric product of $X$
obtained as the quotient of $X^n$ by the permutation action of
$\mathfrak{S}_n$ on factors. Here $\sp{0}X$ reduces to basepoint and
$\sp{1}X=X$. There is a topological embedding
$\sp{n-1}X\hookrightarrow\sp{n}X$ which adjoins the basepoint to a
configuration in $\sp{n-1}X$ and we write $\bsp{n}X:=\sp{n}X/\sp{n-1}X$ for
the cofiber of this embedding (also known as the \textit{reduced} symmetric
product). It is easily seen that $\bsp{n}X$ is the symmetric smash product
$X^{(n)}/\mathfrak{S}_n$, where $X^{(n)}$ is the $n$-fold smash product of $X$
with itself (see \S\ref{constructions}).

The next key result describes $\sj{n}(X)$ and thus ${\mathcal B}_n(X)$ completely
after one unreduced suspension.

\bth\label{main2} There is a homeomorphism\  $\Sigma \sj{n}(X) =
\bsp{n}(\Sigma X)$.
\end{theorem}

Symmetric products being well understood constructions, Theorem \ref{main2}
gives a fairly complete understanding of the stable homotopy type of both
$\sj{n}(X)$ and ${\mathcal B}_n(X)$ and allows for extensive homology
computations. In particular, classical considerations show that $H_*({\mathcal
  B}_n(X))$ is a \textit{direct summand} in
  $$H_*({\mathcal B}_n(X))\hookrightarrow \bigotimes H_{*+1}(K(\tilde
  H_i(X),i+1 ))$$
  where $K(G,i)$ is the Eilenberg-MacLane space whose only non-trivial
  homotopy group is $G$ in degree $i$ (see appendix).  Based on this, we
  indicate in \S\ref{manifolds} for example how to recover some of the
  homological calculations of \cite{bahri} alluded to earlier.  Note that
  combining Theorem \ref{main2} with Theorem \ref{cor1} produces sharp
  connectivity bounds for the reduced symmetric products of a simplicial
  complex (see \cite{braids}).

A few useful corollaries are stated next. We
recall that $X$ is a path-connected space throughout.
We denote by $X^{\vee n}$ the one point union of $n$ copies of $X$.

\bco\label{consequences}\ \\
(a) (Corollary \ref{echar})
If $\chi (X)$ be the Euler characteristic of $X$, then
$$\chi ({\mathcal B}_n(X)) = 1-{1\over n!}(1-\chi (X))
\cdots (n-\chi (X))$$
This has been obtained in the case of topological surfaces in \cite{malchiodi}.\\
(b) (Corollary \ref{circle} )
$\sj{n}(S^1)$ is homeomorphic to $S^{2n-1}\ (\footnote{
According to \cite{vassiliev}, corollary \ref{consequences},
(a) is a result of Caratheodory. No further explicit reference was given.})$.\\
(c) (Lemma \ref{closedsur})
If $C_g$ denotes a closed Riemann surface of genus $g\geq 0$, then
$$\Sigma {\mathcal B}_2(C_g)\simeq
(S^4)^{\vee (2g^2+g)}\vee (S^5)^{\vee 2g}\vee\Sigma^4\bbr P^2$$
(d) (Lemma \ref{product})
There is as well a stable splitting of ${\mathcal B}_2(X\times Y)$
into a six term bouquet.\\
\eco

Our last and one of our most interesting results is that it is possible to
give a streamlined generalization of the James-Thomas-Toda-Whitehead
equivalence (\ref{jmea}) to higher dimensional spheres. Key constructions from
\cite{su} happen to be tailor made for such a generalization.

\bth\label{barofspheres} There is a homotopy equivalence
${\mathcal B}_n(S^k)\simeq\Sigma^{k+1}Q_{n,k}$ where $Q_{n,k}$ is the
quotient of $S^{(k+1)(n-1)-1}$ the unit sphere in the linear
subspace $\{(v_1,\ldots, v_n)\in (\bbr^{k+1})^n\ |\ \sum v_i = 0\}$
under the $\mathfrak{S}_n$-action given by permuting the $v_i$.
\end{theorem}

This theorem is established in \S\ref{spheres}, and in an appendix
we completely determine the
homology of the barycenter spaces of the $2$-sphere.

We conclude this introduction by applying Theorem \ref{barofspheres} to
give a short novel proof of (\ref{jmea}).

\bco\label{maincor}
\cite{james} $\sj{2}(S^{k}) = \Sigma^{k+1}\bbr P^{k}$. \eco

\noindent\proof
When $n=2$ in Theorem \ref{barofspheres}, $Q_{2,k}$ is the unit sphere $S^k$
in $W=\{(v,-v)\in (\bbr^{k+1})^2\}$. The generator of $\bbz_2$ acts on $W$ by
permuting $v$ and $-v$ and hence is multiplication by $-1$ on that
sphere. This is the antipodal action and the claim is immediate.
\hfill$\Box$\vskip 5pt

\vskip 5pt {\sc Acknowledgement:} We would like to thank Abbas Bahri for
suggesting this beautiful problem. We also thank Francesca de Marchis and
Paolo Salvatore for remarks regarding the appendix and \S\ref{barycenters}
respectively.

%***********************************************************************

\section{Notations and Definitions}\label{constructions}

This section sets up notation and basic constructions.  Our spaces are assumed
to be locally compact, hausdorff and unless otherwise stated connected (see \cite{cohen}).
For based spaces $X$ and $Y$ with basepoints $*$,
we write $X\vee Y$ the ``wedge" of $X$ and $Y$ consisting of
$\{(x,y)\in X\times Y\ | \ x=*\ \hbox{or}\ y=*\}$.
We write $X\wedge Y$ the cofiber of the inclusion
$X\vee Y\hookrightarrow X\times Y$.
The ``smash" product $X\wedge Y$ is naturally a based
construction. Throughout
$(D^n,S^{n-1})$ refers to the pair (unit closed disk in $\bbr^n$, its boundary
sphere) and $I = D^1 = [0,1]$.   We write the standard $(n-1)$-simplex as
\begin{equation}\label{simplex}
\Delta_{n-1} = \{ (t_1,\ldots, t_n), t_i\in [0,1], \sum t_i=1\}
\end{equation}
The faces of this simplex correspond to when $t_i=0$. Of course
$\Delta_0=\{1\}$.

\vskip 10pt
\noindent{\sc The Iterated Join.}\\
Given two connected topological spaces
$X$ and $Y$, the join of $X$ and $Y$ is the space
of all segments ``joining points" in $X$ to points in $Y$. It is
denoted by $X*  Y$ and is the identification space
\begin{equation}\label{joindef}
X*  Y := X\times I\times Y/ (x,0,y)\sim (x',0, y), \
(x,1,y)\sim(x,1,y')\ \ \ \ \ \forall x,x'\in X,\ \forall y,y'\in Y
\end{equation}
The homotopy type of $X*  Y$ depends only on the homotopy type of
$X$ and $Y$ (see \cite{brown}). Notice that the join of two based
spaces $(X,x_0)$ and $(Y,y_0)$ is naturally based at $[x_0,{1\over
2}, y_0]$.

The join construction $* $ can be iterated and since it is associative we can
write $X_1* \ldots\cdots* X_n$ the result of performing $*$ $n$-times. This is
the same as the quotient construction
\begin{equation}\label{iterate1}
[X_1\times\cdots\times X_n\times\Delta_{n-1}]_{/\sim}
\end{equation}
where
$$(x_1,\ldots,x_i,\ldots, x_n; t_1,\ldots,t_i= 0, t_n)\sim
(x_1,\ldots, x'_i,\ldots x_n; t_1,\ldots, t_i=0,\ldots, t_n)$$ for
all $x_i,x'_i\in X_i$. When $X_1=\cdots = X_n=X$ we write $X^{* n}$
for the space so obtained.  It is handy to write elements of $X^{*
n}$ formally as $[t_1x_1, t_2x_2,\ldots, t_nx_n]$ with the
understanding that $\sum t_i=1$, $t_i\in [0,1]$ and that $0x = 0x'$
for all $x,x'\in X$.
Note that when collapsing the subspace of points of the form $[t_1x_0,\ldots,
t_nx_0]$ from $X^{*n}$, we obtain the \textit{reduced join}, a space homotopy
equivalent to $X^{* n}$ (see \cite{vogt}).

The cone on $X$ is
$$CX := I\times X/(0 ,x)\sim (0 ,x')\ =\ X*\{v\}$$
viewed as the join of $X$ with a disjoint \textit{vertex} $v$. We will write elements of $CX$
as $[t,x]$ with $[0,x]=[0,x']$. The base $X\hookrightarrow CX$ corresponds
to when $t=1$ and the vertex corresponds to when $t=0$.

The \textit{unreduced} suspension is defined to be $\Sigma X = S^0 *
X$ or equivalently is the quotient $CX/X$ obtained from collapsing
out the base of the cone. Iterating this suspension $n$-times yields
$\Sigma^nX=S^n*X$. The \textit{reduced} suspension on the other
hand is the smash product $S^1\wedge X$.
Since a few of our statements involve homeomorphism
type, we must distinguish between both constructions.

\bex\label{examples} The following relevant facts and alternative descriptions
will be useful:
\begin{enumerate}
\item The join is the strict pushout $X*Y \cong CX\times Y\cup_{X\times Y}X\times CY$.
\item There is a homeomorphism $X*Y \cong CX\times Y/_\sim$\ , where
$(x,y)\sim (x,y'), \forall y,y'\in Y$ and $\forall x\in X$ the base of the cone.
This implies the homotopy equivalence $X*Y\simeq\Sigma X\wedge Y$.
\item  $X^{*n-1}$ embeds in $X^{* n}$ via the map
$[t_1x_1,\ldots, t_{n-1}x_{n-1}]\mapsto [t_1x_1,\ldots,
t_{n-1}x_{n-1},0x_0]$. The image of this embedding includes in the
cone $x_0* X^{*n-1}$ and hence is contractible in $X^{*n}$.
\item If $C$ is
a closed subset of $\bbr^m$ then $\bbr^{m+n}-C$ is homotopy
equivalent to the join $S^{n-1}*(\bbr^m - C)$.%(see \cite{igusa}, lemma 3.1).
\end{enumerate}
\eex

\vskip 5pt\noindent{\sc Further Notation.}\ \\
We write $X^{(n)}:=X\wedge\cdots\wedge X$ the $n$-fold smash product
obtained from $X^n$ after collapsing out the subspace of all tuples
containing the basepoint. This is a space with a ``canonical''
basepoint. We write the
\textit{fat diagonal} in $X^{n}$ (similarly in $X^{(n)}$) as
\begin{equation}\label{fatdiagonal}
\Delta_{fat}:=\{[x_1,\ldots, x_n ]\in X^{n}\ |\ x_i=x_j\ \hbox{for
some}\ i\neq j\}
\end{equation}

A \textit{based} action of a group $G$ on a pointed space $Y$ is an
action that fixes the basepoint. The action is based free if it is
free away from that basepoint.  If $X$ admits a based right action
of $G$ and $Y$ a based left action, then we can take the diagonal
quotients $X\times_GY$, $X\vee_GY = (X/G)\vee (Y/G)$ and
$X\wedge_GY$ (the quotient of the previous two).  The prototypical
example of a based action we consider in this paper is that of the
symmetric group $G=\mathfrak{S}_n$ acting on the $n$-fold smash
$Y=X^{(n)}$ by permuting factors (the basepoint is the canonical
basepoint).

The ``homotopy colimit or pushout" of two maps $f: Y\lrar X, g :
Y\lrar Z$ is the double mapping cylinder
$$X\sqcup (Y\times [0,1])\sqcup Z\ /_{\sim}\ \ ,\ \
(y,0)\sim f(y)\ ,\ (y,1)\sim g(y),\ \forall y\in Y
$$
This is written $\hbox{hocolim} (X\fract{f}{\la 2} Y\fract{g}{\lrar} Z)$ so
that for instance, Example \ref{examples}, (1), becomes
$$
X*Y = \hbox{hocolim}(X\leftarrow X\times Y\rightarrow Y)
$$

%********************************************************

\section{Symmetric Joins of Spaces}\label{symjoin}

The symmetric group on $n$-letters $\mathfrak{S}_n$ acts on $X^{*
n}$ by permuting factors
$$\sigma [t_1x_1,\ldots, t_nx_n] =
[t_{\sigma (1)}x_{\sigma (1)},\ldots, t_{\sigma (n)}x_{\sigma (n)}]\
\ ,\ \ \sigma\in\mathfrak{S}_n
$$
The quotient under this action is the {\sl $n$-th symmetric join}

$\sj{n}(X) := X^{*  n}/_{\mathfrak{S}_n}$.
Alternatively and based on (\ref{iterate1}) we can write
\begin{equation}\label{identification1}
\sj{n}(X) := \coprod_{k=1}^n
\Delta_{k-1}\times_{\mathfrak{S}_k}X^k/_{\sim}
\end{equation}
where $\mathfrak{S}_k$ acts on $X^k$ and on
$\Delta_{k-1}\subset\bbr^k$ by permuting factors, and
$\times_{\mathfrak{S}_k}$ means taking the quotient with respect to the
diagonal $\mathfrak{S}_k$ action. The equivalence relation $\sim$ indicates
further identifications of the form
$$\hbox{(i)}\ \ \ \ \ \ \ \
[x_1,\ldots,x_k; t_1,\ldots, t_k]\sim [x_1,\ldots, \hat x_i,\ldots x_k;
t_1,\ldots, \hat t_i,\ldots, t_k]$$
whenever $t_i=0$, here \textit{hat} means deletion.
Note that these identifications are made along the
boundary of the simplex $\Delta_{k-1}$.

\vskip 5pt \noindent{\sc Notation}:
Elements of $\sp{n}X$ will be written as \textit{formal abelian} sums $\sum x_i =
x_1+\cdots +x_n$ and this represents the orbit of $(x_1,\ldots, x_n)\in X^n$
under the $\mathfrak{S}_n$ action.  Similarly elements
$\zeta$ in $\sj{n}(X)$ are written in the form $\sum
t_ix_i$, $\sum t_i= 1$, with the understanding that if one of the
$t_i$'s is $0$, then that entry is \textit{suppressed} from the sum.
This writing is unique if none of the $t_i$'s is zero.

\bre By construction we have an inclusion
$\sj{n-1}X\hookrightarrow\sj{n}X$. This subspace is contractible
in $\sj{n}X$. If $x_0$ is a basepoint in $X$, it is fairly transparent
how to write the contraction
\begin{eqnarray*}
I\times \sj{n-1}(X)&\lrar& \sj{n}(X)\\
(t, \sum t_ix_i) &\longmapsto&tx_0 + \sum (1-t)t_ix_i
\end{eqnarray*}
In particular the infinite symmetric join is weakly contractible (i.e.
trivial homotopy groups).\ere

\bex\label{jointocone} There is an embedding
$\sj{n}(X)\hookrightarrow \sp{n}(CX)$ which sends $\sum_1^n t_ix_i$
to $\sum [t_i,x_i]$. This map is well-defined because the vertex of
the cone is at $t=0$. The space $\sp{n}(CX)$ is of course
contractible since $CX$ is. In the case of the circle $X=S^1$ and
$CS^1=D^2$, we obtain an explicit embedding of $\sj{n}(S^1)$ into
Euclidean space $\sp{n}(D^2) \cong D^{2n}$. \eex

It becomes now an interesting problem to determine the homotopy type
of symmetric joins of some standard spaces.  An early cute
characterization was given in \cite{james} for the ``symmetric square" of spheres

\ble\label{jmea1}\cite{james} There is a homeomorphism $\sj{2}(S^n ) =
\Sigma^{n+1}\bbr P^n$.
\ele

\noindent\proof  This proof is so short and clever that we reproduce it.
Let $D=D^{n+1}$ be the unit $(n+1)$-disc and
identify it with $D^n=CS^{n-1}$ with vertex at the origin. By
inserting $D^n = CS^{n-1}=S^{n-1}\times I/\sim$ in the definition of
the join (\ref{iterate1}), we can rewrite $S^{n-1}* X$ as $D^n\times
X/(x,y)\sim (x,y'),\ x\in S^{n-1}$ (that is $\theta\times
X\sim\theta$ for $\theta\in\partial D^n$).

Write $(S^n)^{* 2}$ as $S^n\times D^{n+1}\cup D^{n+1}\times S^n$ as in
(Example \ref{examples}, (1)). To $(x,y)\in D\times\bbr P^n$ draw a line through $x$
parallel to $y$. This meets $S^{n}=\partial D$ in points
$z,z'$. Choose $z$ to be closest to $x$, and choose $w$ (between $x$ and $z'$)
such that $x$ is the midpoint of $z$ and $w$. We define
$$f : D\times \bbr P^n\lrar \sj{2}(S^n) \ \ ,\ \ (x,y)\mapsto [(z,w)]$$
If $z$ and $z'$ are equidistant to $x$, $f(x,y)=[(z,z')]$ and there is
no ambiguity so our map is well-defined.
Notice that if $x\in S^n$, then $z=w=x$ (is independent of $y$). We can
then factorize the map $f$ through the quotient $D\times \bbr P^n/(x,y)\sim
(x,y'),\forall x\in \partial D$.  This quotient is another description
of $\Sigma^{n}* \bbr P^n = \Sigma^{n+1}\bbr P^n$ (see Example \ref{examples},
(2)).  The map so obtained
$D\times\bbr P^n/_{\sim}\rightarrow\sj{2}(S^n)$ is continuous and
bijective by construction, hence is bicontinuous (being closed).
\hfill$\Box$\vskip 5pt

Note in particular that $\sj{2}(S^1)\cong S^3$ which is a special case of corollary \ref{circle}. For more general spaces $X$ one has the following global description
of $\sj{2}(X)$.

\ble\label{square} There is a homotopy pushout
$$\sj{2}X = \hbox{hocolim} (X\fract{p_2}{\la 2} X^2\fract{\pi}{\lrar}\sp{2}X)$$
where $p_2$ is the projection onto
the second factor, and $\pi$ is the $\bbz_2$-quotient map.  \ele

\noindent\proof  An element of $\sj{2}(X)$ is written as $t_1x + t_2y$,
$t_1+t_2=1$ with the identification $0x+1y =y$. By using the order on
the $t_i$'s in $I$, we can write
\begin{eqnarray*}
\sj{2}(X) &=& \{(t_1,t_2,x_1,x_2)\ | \ t_1\leq t_2, t_1+t_2=1\}/_\sim\\
&=& I\times (X\times X)/_\sim
\end{eqnarray*}
where $I=\{0\leq t_1\leq t_2\leq 1, t_1+t_2=1\}$ is a copy of the
one-simplex, and the identification $\sim$ is such that $(0,1,x,y)\sim
y$ and $({1\over 2},{1\over 2}, x,y)\sim ({1\over 2},{1\over 2},
y,x)$. But $(0,1)$ and $({1\over 2},{1\over 2})$ are precisely the
faces or endpoints of $I$ and hence the claim.
\hfill$\Box$\vskip 5pt

Pictorially this double mapping cylinder can be depicted as in the figure
$$\xymatrix{X\ar@{-}[rrr]^{\hbox{X$\times$ X}}&&&\sp{2}X}$$
As is clear it is possible to give a multiple pushout
description of $\sj{n}(X)$ for $n\geq 3$ as in lemma \ref{square}.
For the case $n=3$ this proceeds as follows. Write elements as
$(t_1,t_2,t_3;x_1,x_2,x_3)$ with $t_1+t_2+t_3=1$, $t_1\leq t_2\leq
t_3$.  The faces of this 2-simplex are: $t_1=0$ ( and hence
$t_2+t_3=1$ forming the ``right edge''), $t_1=t_2$ and ${1\over
3}\leq t_3\leq 1$ (the ``left edge''), and thirdly $t_2=t_3$ and
$0\leq t_1\leq {1\over 3}$ (bottom edge). The top vertex corresponds
to when $t_1=t_2=0, t_3=1$, the bottom left as when
$t_1=t_2=t_3={1\over 3}$ and the bottom right vertex to when $t_1=0,
t_2=t_3={1\over 2}$. The multiple homotopy pushout diagram representing $\sj{3}(X)$
is depicted below.
$$\small\xymatrix{&&X\ar@{-}[ddrr]^{\hbox{X$\times$ X}}\\
&&X\times X\times X\\
\sp{3}X\ar@{-}^{\sp{2}\hbox{X$\times$ X}}[uurr]\ar@{-}[rrrr]^{\hbox{X$\times\sp{2}$X}}_{{\phantom{\vdots}
\atop\hbox{{\bf Figure 1}: the colimit diagram for $\sj{3}(X)$}}}&&&&\sp{2}X}$$

The construction for $n>3$ is predictable. Such a description is
what we need to establish the simple connectivity of symmetric joins
which is a first step towards Theorem \ref{cor1}.

\ble\label{specialpushout} The following homotopy pushout is simply
connected where the gluing on the left is via projection $p_1$ on
the first coordinate, and on the right via the concatenation product $\mu$
$$\xymatrix{
\hbox{SP}^i\hbox{X}\ar@{-}[rrr]^{\sp{i}\hbox{X}\times\sp{j}\hbox{X}\ \ \ \ \ \ }&&&\sp{i+j}X
\ \ \ \ \ \ i,j\geq 1}$$
\ele

\noindent\proof  By the Van-Kampen theorem, the fundamental group $G$ of the pushout
has generators from $\pi_1(\sp{i}X)$ and $\pi_1(\sp{i+j}X)$ subject to the
``VK-relations" demanding in particular
that $\mu_*(x)$ and $p_{1*}(x)$ map to the same element
in $G$ for every $x\in\pi_1(\sp{i}X\times\sp{j}X)$.
It is an easily established fact \cite{liao} that the connectivity of
  symmetric products is at least the same as that of the underlying space.  If
  $X$ is simply connected then so are all the spaces involved along the edge
  and endpoints of the pushout diagram and hence the pushout is also simply
  connected. In the general case one can use another useful property of
  symmetric products which is that $\pi_1(\sp{i}(X)) \cong H_1(X;\bbz )$ for $i\geq 2$,
  so that in particular the basepoint inclusion induces a map
  $\pi_1(X)\lrar\pi_1(\sp{i}X)$ which is surjective for $i\geq 2$.
  It follows from there that any element of $\pi_1(\sp{i+j}(X))$ is
  in the image of an element from $\pi_1(\sp{j}(X))$ under the
  concatenation product
  $\mu : \sp{i}(X)\times\sp{j}(X)\lrar\sp{i+j}(X)$. But this element projects trivially
in $\pi_1(\sp{i}X)$ under $p_{1*}$ so that its image in the pushout $G$ in trivial by the VK relations. On the other hand, these same relations show that any
  generator from $\pi_1(\sp{i}X)$ is identified to a generator coming from $\pi_1(\sp{i+j}X)$
  so that its image in $G$ is also trivial. This proves our claim.
\hfill$\Box$\vskip 5pt

\bth\label{simplyconnected} $\sj{n}(X)$ is simply connected for
$n\geq 2$.\end{theorem}

\noindent\proof
  According to lemmas \ref{square} and \ref{specialpushout}, $\sj{2}(X)$
  is simply connected.  We can then proceed by induction.  Consider
  $\sj{3}(X)$ as in figure 1.  Label the vertices $V_1:=X$, $V_2=\sp{2}(X)$,
  $V_3=\sp{3}(X)$ and their opposite edges $E_1,E_2,E_3$, so for instance
  $E_1 = X\times\sp{2}(X)$.
 Let $U_i = \sj{3}(X)-V_i$. The $U_i$'s are open subspaces and
  each $U_i$ deformation retracts onto $E_i$ (this is proved further
  down). But each $E_i$ is a space obtained by iterative pushout
  and is simply connected by lemma
  \ref{specialpushout} and our inductive hypothesis.
  Since $\sj{3}(X)=U_1\cup U_2\cup U_3$ is covered by
  three open simply connected subsets with non-empty intersection, Van-Kampen
  theorem implies that $\sj{3}(X)$ is simply connected.

  The proof for $n=3$ extends with little change to larger $n$.  The pushout
  diagram for $\sj{n}(X)$ is a labeled $n-1$-simplex with vertices
  $V_i=\sp{i}(X)$, $1\leq i\leq n$. Opposite to each $V_i$ is a face
  representing a pushout $E_i$. The open set $U_i=\sj{n}(X)-V_i$ deformation
  retracts onto $E_i$ and so it is simply
  connected by induction. Since the $U_i$ cover $\sj{n}(X)$, the claim follows
  by Van-Kampen.

  It remains to establish that indeed $U_i$ deformation retracts onto
  $E_i$. Figure 1 suggests an argument of proof.  Pick a vertex $v_i$ and write
  $\Delta_n=\Delta_{n-1}*\{v_i\}$ a join, then
  there is an obvious geometric deformation of $\Delta_n-\{v_i\}$ onto
  $\Delta_{n-1}$ given by $r_i: (\Delta_n-\{v_i\}) \times
  I\lrar \Delta_n-\{v_i\}$,
$$r_i(sx + (1-s)v_i,t) = (1-t)(sx+(1-s)v_i)+ tx$$
with $x\in \Delta_{n-1}$, $sx+(1-s)v_i\in \Delta_{n-1}*v_i$ and
$t\in I$. This homotopy is well-defined since $s\neq 0$ (and hence
$x$ can be chosen) and it retracts the complement of the vertex
$v_i$ onto its opposite edge leaving the other faces invariant. Let
$U_i$ be the image of
$$\phi_i : (\Delta_n-\{v_i\})\times X^{n+1}\lrar {\mathcal B}_{n+1}(X)$$
Then post composing $\phi_i$ with the map
\begin{eqnarray*}
(\Delta_n-\{v_i\})\times X^{n+1}\times
I&\lrar&(\Delta_n-\{v_i\})\times
X^{n+1}
\end{eqnarray*}
which is $r_i$ on $\Delta_n-\{0\}\times I$ and the identity on $X^{n+1}$,
gives the desired  retraction of $U_i$ onto $E_i$.
\hfill$\Box$\vskip 5pt

%********************************************************

\section{Symmetric Joins and Symmetric Products}

The multiple pushout description of $\sj{n}(X)$, although useful for the
fundamental group calculation, doesn't yield itself to an easy homological
analysis. A much finer and surprisingly simple description is given next.

\bpr\label{key}
For spaces $X_1,\ldots, X_n$, there is a homeomorphism
$$C(X_1*  X_2*  \cdots *  X_n)\cong CX_1\times\cdots\times CX_n$$
In particular $C(X^{*  n})\cong (CX)^n$. This homeomorphism is
$\mathfrak{S}_n$-equivariant so that
$$C(\sj{n}(X))\cong \sp{n}(CX)\ \ \ \hbox{and}\ \ \
\Sigma \sj{n}(X)\cong\bsp{n}(\Sigma X)$$
\epr

\noindent\proof  Start by writing elements of $CX$ as $[t,x]$
  with $[0,x]=[0,x']$.  The base of the cone corresponds to $X\hookrightarrow
  CX$, $x\mapsto [1,x]$.  Let's view $X_1* \cdots*
  X_n$ as the quotient of $X_1\times\cdots\times X_n\times\Delta_{n-1}$ as in (\ref{iterate1}),
  consider the map
\begin{eqnarray*}
\Phi : C(X_1*  X_2* \cdots*  X_n)&{\ra 3}& \prod_{j=1}^nCX_j\\
\ [s, (x_1,\ldots, x_n; t_1,\ldots, t_n)]&\lrar&
([s{t_1\over t}, x_1],\cdots
[s{t_n\over t}, x_n])
\end{eqnarray*}
where $t = max_j\{t_j\}$. Notice that $t\neq 0$ since $\sum t_i=1$
and $\Phi$ is well-defined. This map is a homeomorphism with inverse
$$\Phi^{-1}([s_1,x_1],\ldots, [s_n,x_n])
= \begin{cases}[max_j\{s_j\}, (x_1,\ldots, x_n; {s_1\over s},\ldots, {s_n\over s})],&
\ \hbox{if}\ max_j\{s_j\} > 0 \\
*,&\ \hbox{if}\ max_j\{s_j\}=0
\end{cases}$$
where we have written $s = s_1+\cdots +s_n$. This map is easily seen to be
continuous.

Note that $\Phi$ carries the base of the cone on $X_1* X_2* \ldots* X_n$
homeomorphically onto the subspace
$$D_n := \{([t_1,x_1],\ldots, [t_n,x_n]\ |\ t_j=1\ \hbox{for at least one $j$}\}\subset
\prod_{j=1}^n CX_i$$
This is saying that there is a homemorphism of pairs
$$\Phi : (C(X_1*  X_2* \ldots*  X_n ),X_1*  X_2* \ldots*  X_n)
\cong (\prod_j CX_j, D_n)$$
and so collapsing out subspaces we obtain a homeomorphism
$$\Phi : \Sigma (X_1*  X_2* \ldots*  X_n )\cong \Sigma X_1\wedge\cdots\wedge
\Sigma X_n$$ where $\Sigma X = CX/X$. When all the $X_j$'s equal to $X$, the
map $\Phi$ is evidently $\mathfrak{S}_n$-equivariant and there is an induced
homeomorphism $\Sigma (\sj{n}X)\cong\bsp{n}(\Sigma X)$.
\hfill$\Box$\vskip 5pt

We give an alternate proof of the homotopy equivalence
$\bsp{n}(\Sigma X)\simeq\Sigma \sj{n}(X)$
in Corollary \ref{suspendaction}.

\bre The earliest occurence of the homeomorphism
$\Sigma (X*X)\cong \Sigma X\wedge\Sigma X$ that we are aware of
is in \cite{cohen}. The above proof is similar to Cohen's argument.
The homeomorphism $C(X)\times C(Y) = C(X* Y)$
is quite standard and can be seen geometrically by thinking of
$CX\times CY$ as a square with center the vertex of a cone on the
the boundary of this square $CX\times Y\cup X\times CY$. \ere

As an immediate corollary of proposition \ref{key} and of
(\ref{jmea}) we obtain the identification
$\bsp{2}(S^n) = \Sigma^{n+1}\bbr P^{n-1}$.
For a generalization see Proposition \ref{su}.
Another immediate corollary is the following.

\bco\label{circle} $\sj{n}(S^1)\cong S^{2n-1}$.  \eco

\noindent\proof
  The cone $CS^1$ is homeomorphic to the two dimensional closed disc $D^2$ so
  that $\sp{n}(D^2)\cong C(\sj{n}(S^1))$. It is well-known that
  $\sp{n}D^2\cong D^{2n}$ with boundary $\partial\sp{n}(D^2) = \sp{n}D - \sp{n}(int(D))$
  the sphere $S^{2n-1}=\partial D^{2n}$ as is checked in \cite{ks}. This
  boundary sphere must correspond to unordered tuples with one entry in $\partial D$.
Under the homeomorphism constructed in Proposition \ref{key},
 this is precisely the base of the
 cone on $\sj{n}(S^1)$ and the claim follows.
\hfill$\Box$\vskip 5pt

Another fairly useful corollary is the calculation of the Euler characteristic of
symmetric joins (and hence of barycenter spaces in light of Theorem
\ref{main1}).

\bco\label{echar} If $\chi (X)$ is the Euler characteristic of a connected
space $X$, then
$$\chi (\sj{n}(X)) = 1-{1\over n!}(1-\chi (X))
\cdots (n-\chi (X))$$
\eco

\vskip 5pt\noindent\proof We recall that
$\chi (\Sigma X) = 2-\chi (X)$, and that\\
(i) $\chi (\sp{n}X) =
{1\over n!}(\chi (X)+n-1)\cdots (\chi(X)+1)\chi (X))$. This is extracted from
the formula of MacDonald \cite{macdonald} that the Euler characteristics for
all $n$ put together form the series
$$\sum_{n\geq 0}^\infty \chi (\sp{n}X)t^n = \left({1\over 1-t}\right)^{\chi (X)}$$
so that $\chi (\sp{n}X)$ is the number of non-negative solutions of
$j_1+j_2+\cdots +j_{\chi (X)} = n$; that is ${\chi (X)+n-1\choose \chi (X)-1}$.
(ii) $\chi (A\cup B) = \chi (A)+\chi (B)-\chi (A\cap B)$. In particular
if $A\hookrightarrow X$ is a cofibration, then the quotient
$X/A$ is identified with the mapping cone $X\cup_ACA$, and
hence $\chi (X/A) = \chi (X) + 1 - \chi (A)$. Applying this formula to
$\bsp{n}X = \sp{n}X/\sp{n-1}X$, we see that
$$\chi (\bsp{n}X) = 1+ {\chi (X)+n-1\choose \chi (X)-1} -
{\chi (X)+n-2\choose \chi (X)-1}$$
From this we obtain that
$$\chi (\sj{n}(X)) =
2-\chi (\bsp{n}\Sigma X) = 1 - {n+1-\chi (X)\choose 1-\chi (X)} +
{n-\chi (X)\choose 1-\chi (X)}$$ and this gives precisely our formula.
\hfill$\Box$

\bre It is not surprising that $\bsp{n}(\Sigma X)\cong\Sigma\sj{n}(X)$
is again a suspension. Indeed
let $Y$ be a co-$H$ space with comultiplication $\nabla : Y\lrar Y\vee
Y$ (by definition of co-$H$ space, this means that the composite
$Y\fract{\nabla}{\lrar}Y\vee Y\hookrightarrow Y\times Y$ is homotopic
to the diagonal map). Then we claim that $\bsp{n}(Y)$ is also co-H.
To see this set $\bsp{0}(Y)=S^0$. It is not hard to check that
\begin{equation}\label{spwedge}
\bsp{n}(X\vee Y) = \bigvee_{r+s=n}\bsp{r}(X)\wedge\bsp{s}(Y)
\end{equation}
The co-$H$ space structure on $\bsp{n}(Y)$ is obtained from the composite
$$\bsp{n}(Y)\lrar\bsp{n}(Y\vee Y)=
\bigvee_0^n\bsp{i}(Y)\wedge\bsp{n-i}(Y) \lrar \bsp{n}Y\vee\bsp{n}Y
$$ with the second map being projection on the wedge summands corresponding
to when $i=n$ and $i=0$.
\ere

%***************************************************************

\section{Barycenter Spaces}\label{barycenters}

We now turn to the spaces of most interest in this paper. The $n$-th
barycenter space ${\mathcal B}_n(X)$ is defined as a quotient of the symmetric
join $\sj{n}(X)$ under one further type of identifications of the form
\begin{equation}\label{pointaddition}
\hbox{(ii)}\ \ \ \ \ \ \ t_1x_1+t_2x_1+t_3x_3+\cdots +t_nx_n\sim
(t_1+t_2)x_1+t_3x_3+\cdots +t_nx_n
\end{equation}
More precisely we can write
\begin{equation}\label{identification2}
{\mathcal B}_{n}(X) := \coprod_{k=1}^n
\Delta_{k-1}\times_{\mathfrak{S}_k}X^k/_{\approx}
\end{equation}
where $\approx$ consists of the identifications (i) as in
(\ref{identification1}) and (ii) as in (\ref{pointaddition}).

Here too, it is useful to think of a
point of ${\mathcal B}_n(X)$ as an abelian sum $t_1x_1+ \cdots + t_nx_n$,
$\sum t_i=1$, with
the topology that when $t_i=0$ the entry $0\cdot x_i$ can be discarded
from the sum, and when $x_i$ moves in coincidence with $x_j$,
the coefficients add up. As before, there is an
embedding ${\mathcal B}_{n-1}(X)\lrar {\mathcal B}_n(X)$ with  contractible
image.

\bre If $X$ is embedded in some Euclidean space $\bbr^n$, then an element
$\sum t_ix_i$, $x_i\neq x_j$ can be viewed as the ``barycenter" of
the $x_i$'s with corresponding weights $t_i$ and the resulting space
of all weighted barycenters is like \textit{a convex hull} in
$\bbr^n$. Ensuring that geometrically these weighted barycenters
are distinct for different choices of finite tuples $(x_1,\ldots,
x_n)$, is requiring that $X$ is embedded in
``general position" in $\bbr^{N}$ for large enough $N$.  For example
${\mathcal B}_2(\bbr^1)$ can be identified with the
convex hull of the Veronese embedding
$\bbr^1\hookrightarrow\bbr^{N}$, $t\mapsto (t,t^2,t^3,\ldots ,
t^N)$ for $N\geq 3$.
\ere

\bex (Bunch of points) Let ${\bf n}$ be a finite set
which we can assume embedded in Euclidean space as the
vertices of an $n-1$-simplex $\bf\Delta$. Then
${\mathcal B}_k({\bf n})$ is the union of all
$(k-1)$-faces of $\bf\Delta$. This is
a bouquet of spheres $\bigvee^{} S^{k-1}$. For example
${\mathcal B}_{n-1}({\bf n})= S^{n-2}$.  \eex

In this section we give as simple as possible of a description of
the homotopy type of ${\mathcal B}_n(X)$ and observe that it is
equivalent to $\sj{n}(X)$.
Recall again that $X^{(n)}$ is the quotient of $X^n$ by the fat wedge
of all tuples with at least one entry at $x_0$.

\bth\label{equivalence} For a based connected simplicial complex $X$,
there is a homotopy equivalence
$${\mathcal B}_n(X)\simeq
 S^{n-1}\wedge_{\mathfrak{S}_n}X^{(n)}$$
 with $\mathfrak S_n$ acting on $S^{n-1}\subset\bbr^n$ by permuting
 coordinates, fixing the point ${1\over\sqrt{n}}(1,\ldots, 1)$.
\end{theorem}

\noindent\proof  With $\mathcal B_n(X)$ as in (\ref{identification2}),
and $x_0\in X$ the basepoint, consider the subspace
$$W_n := \left\{\sum t_ix_i\in {\mathcal B}_{n}(X), x_i=x_0\ \hbox{for some}\ i\right\}$$
Notice that ${\mathcal B}_{n-1}X$ includes in $W_n$ as the subspace of all points of the form $\sum_1^{n-1} t_ix_i + 0x_0$.
We claim that $W_n$ is contractible. We construct a
contraction as follows. Write an element of $W_n$ in the form $\zeta
= t_1y_1+\cdots +t_ky_k + s_{0}x_0$, for some $k< n$, with
 $y_i\neq x_0$. We have a well-defined homotopy $F:
W_n\times I \lrar W_n$
\begin{eqnarray}\label{contraction}
&&(t_1y_1+\cdots +t_ky_k + s_{0}x_0,t)\\
&\mapsto&
(1-t)t_1y_1+\cdots + (1-t)t_ky_k + (s_0+t(1-s_{0}))x_0\nonumber
\end{eqnarray}
This is well-defined because the sum of coefficients is always $1$,
and because it is continuous (we have to check continuity when one of the
$y_i$'s approaches $x_0$ but this is an inspection).
This homotopy is evidently a contraction with $F(0,-) = id$ and $F(1,-)$
the constant function at $x_0$.
Since $W_n$ is contractible, we can
collapse it out without changing the homotopy type.

Now $\approx$ identifies as in (i) the
subspace of elements $\sum t_ix_i$ for which $(t_1,\ldots, t_n)$ lie
in the boundary of $\Delta_{n-1}$ (i.e. for which $t_i=0$ for some
$i$) with a subspace of $W_n\subset
{\mathcal B}_{n}(X)$. Similarly
configurations containing the basepoint $x_0$ are in $W_n$.
What we obtain by identifying $W_n$ to point is then the quotient
\begin{equation}\label{identification3}
{\mathcal B}_{n}(X) \simeq
\left([\Delta\ltimes_{\mathfrak{S}_n}X^{(n)}]/
[\partial\Delta\ltimes_{\mathfrak{S}_n}X^{(n)}]\right)_\approx
= S^{n-1}\wedge_{\mathfrak{S}_n}X^{(n)}/_\approx
\end{equation}
and $\approx$ is the identification (\ref{pointaddition}).
Here $X\ltimes Y := X\times Y/X\times *$ is the ``half-smash" product.
The smash product is
taken with respect to the canonical basepoints in both $X^{(n)}$
and $S^{n-1}:=\Delta/\partial\Delta$. The symmetric group $\mathfrak{S}_n$ acts on $S^{n-1}=\Delta/\partial\Delta$ by
permuting the $t$'s of the simplex.
Now note that (\ref{pointaddition}) identifies
points in the \textit{fat diagonal} $\Delta_{fat}$; that is
the subset of the form $\sum t_ix_i$ with $x_i=x_j$ for some
$i\neq j$, to points in $W_n$ which we have again collapsed out. We can
then write
\begin{eqnarray}\label{identification4}
{\mathcal B}_n(X)&\simeq& (S^{n-1}\wedge_{\mathfrak{S}_n}X^{(n)})/
(S^{n-1}\wedge_{\mathfrak{S}_n}\Delta_{fat})\nonumber\\
&=&  S^{n-1}\wedge_{\mathfrak{S}_n}(X^{(n)}/\Delta_{fat})
\end{eqnarray}
This last equivalence follows from the fact that $S^{n-1}\wedge
X^{(n)}/S^{n-1}\wedge \Delta_{fat} = S^{n-1}\wedge
(X^{(n)}/\Delta_{fat})$ and that this identification passes to the
$\mathfrak{S}_n$-quotients.

It then remains to see that (\ref{identification4}) is the same
as the expression in the theorem.  This
would follow immediately if we knew that
$S^{n-1}\wedge_{\mathfrak{S}_n}\Delta_{fat}$ were contractible. But
this is precisely the content of lemma \ref{contractible} next. To
summarize then, we have shown the string of equivalences
$$\mathcal B_n(X) \simeq \mathcal B_n(X)/W_n\simeq
S^{n-1}\wedge_{\mathfrak{S}_n}(X^{(n)}/\Delta_{fat})\simeq
S^{n-1}\wedge_{\mathfrak{S}_n}X^{(n)}
$$
Finally to get the version of our
theorem, we need replace this sphere and this action by
$S^{n-1}\subset\bbr^n$ and the permutation action on coordinates,
which is possible according to Lemma \ref{deltasphere}.
\hfill$\Box$

\bco\label{equivalence2}
 There is a homotopy equivalence ${\mathcal B}_n(X)\simeq\sj{n}X$.
\eco

We will derive this important corollary in two different ways.

\vskip 5pt\noindent\proof (of corollary \ref{equivalence2})
The first approach relies on a theorem of R. Vogt.
Proposition 1.9 of \cite{vogt} implies the existence of
a $\mathfrak{S}_n$-equivariant equivalence (see remark \ref{vogt})
\begin{equation}\label{vogtresult}
X^{* n}\simeq_{\mathfrak{S}_n} S^{n-1}\wedge X^{(n)}
\end{equation}
Passing to the quotient under the symmetric group action
gives that $\sj{n}(X)\simeq S^{n-1}\wedge_{\mathfrak S_n} X^{(n)}$
and hence the desired result. Here $X^{* n}$ is the reduced join,
 ${1\over\sqrt{n}}(1,\ldots, 1)$ and $(x_0,\ldots ,x_0)$
 are the basepoints in $S^{n-1}\subset\bbr^n$ and $X^n$ respectively.
\hfill$\Box$\vskip 5pt

The proof of the equivalence (\ref{vogtresult}) in \cite{vogt}
takes some space and uses techniques
different from ours. We prefer to give a
second proof of this result that is more self-contained and
close to the proof of Theorem \ref{equivalence}.

\vskip 5pt\noindent\proof (second proof of corollary \ref{equivalence2})
When $n=2$, the proof is immediate since there is a
strict pushout
$$\xymatrix{
I\times X\ar[r]^f\ar[d]^{p_2}&\sj{2}X\ar[d]\\
X\ar[r]&{\mathcal B}_2X}
$$
where $p_2$ is projection onto $X$, and $f (t,x) = tx + (1-t)x$.
Since the left vertical map is an equivalence, then so is the right
vertical one by a standard property of pushouts.

For the more general case $n>2$, we observe that we can follow
exactly the same steps as in the proof of Theorem \ref{equivalence}
for the barycenter spaces and define the subspace
$$V_n:= \left\{\sum t_iy_i\in \sj{n}(X), y_i=x_0\ \hbox{for some}\ i\right\}
$$
then attempt to show that $V_n$ is contractible. The contraction we have
written down for $W_n$ in the proof of Theorem \ref{equivalence} does not however
apply to $V_n$. The problem is that it isn't anymore continuous. To go around this,
we go back to the quotient construction in (\ref{identification1})
and make an extra identification as follows. We keep
the usual identification (i) along the boundaries of the $\Delta_i$
but make the additional identification at the basepoint
$x_0$ of $X$
\begin{equation}\label{x0relation}
tx_0 +sx_0 + ... \sim (t+s)x_0 + ...
\end{equation}
This is of course the identification used in defining $\mathcal B_n$ but we only demand
it holds true at $x_0$. Let's refer to this new
quotient by $T_n(X)$ with projection
$$\pi : \sj{n}X\lrar T_n(X) = \coprod_{k=1}^n
\Delta_{k-1}\times_{\mathfrak{S}_k}X^k/_{\sim}$$
Suppose we show that $\pi$ is a homotopy equivalence. Then we can instead work
with $T_n(X)$ and define the subspace $\pi(V_n)$ which again consists of all
configurations in $T_n(X)$ where $x_0$ appears. We could then show that $\pi(V_n)$ is
contractible since now the contraction
we've written down in (\ref{contraction}) works; i.e. it is continuous
and this is the reason for introducing the relation (\ref{x0relation}).

What is left to do then is show that $\pi$ is an equivalence. To that end
we use a version of the Begle-Vietoris theorem which states that a map between
locally compact, locally contractible spaces such that the inverse image
of any point is contractible, is necessarily a homotopy equivalence
\cite{smale}. Now the preimage of a configuration
$\sum t_iy_i = t_1y_1+\cdots +t_ky_k$
is the same configuration again if no $x_0$ figures
in this sum. In general, the inverse image of
$\zeta = t_1y_1+\cdots +t_ky_k + sx_0$, with $\sum t_i+s=1$ and $y_i\neq x_0$, is
$$\pi^{-1}(\zeta ) = \left\{\sum_1^k t_iy_i + s_1x_0+\cdots +s_{n-k}x_0\ \ \ \ ,\ \ \
s_1+\cdots +s_{n-k} = s= 1 - \sum t_i\right\}$$
This is a copy of a simplex hence contractible. Therefore
$$\sj{n}(X)\simeq T_n(X)\simeq T_n(X)/\pi (V_n)\simeq S^{n-1}\wedge_{\mathfrak S_n}X^{(n)}$$
and the proof is complete.
\hfill$\Box$\vskip 5pt

Finally the following important lemma was needed in the proof of
Theorem \ref{equivalence} above.

\ble\label{contractible} Let $\Delta_{fat}$ be the fat diagonal in
$X^{(n)}$.
Then $S^{n-1}\wedge_{\mathfrak{S}_n}\Delta_{fat}$ is
contractible. \ele

\vskip 5pt\noindent\proof
The justification for this lemma is in Arone-Dwyer \cite{dwyer},
proposition 7.11. We recall their statement. Given a simplicial
$G$-space $Y$, let $Iso(Y)$ denote the collection of all subgroups
of $G$ which appear as isotropy of simplices of $Y$.
Let $G=\mathfrak S_n$ and denote by $\mathcal P$ the subgroups
of $\mathfrak{S}_n$ of the form
$\mathfrak{S}_{n_1}\times \cdots\times\mathfrak{S}_{n_j}$ with $j>1$
and $\sum n_i=n$ with $n_i>1$ for at least one $i$. Let $S^k$ be a
sphere such that $S^k/K$ is contractible for all $K\in{\mathcal P}$.
This is the case for example when $k=n$, $S^n=(S^1)^{(n)}$ and
$\mathfrak{S}_n$ acts by permutation as in the case of \cite{dwyer}.
Then by an elegant
application of equivariant cell attachments, if $Y$ is
a pointed $\mathfrak{S}_n$-simplicial space with $Iso(Y)\subset {\mathcal
P}\cup\{\mathfrak{S}_n\}$, the quotient $S^k\wedge_{\mathfrak{S}_n} Y$
is contractible \cite{dwyer}. We apply this to $Y=\Delta_{fat}\subset X^{(n)}$ and to the sphere $\Delta/\partial\Delta$ with its permutation action of $\mathfrak S_n$.
First we can verify that $\Delta/\partial\Delta$ is $\mathfrak S_n$-equivariantly
homeomorphic to the unit sphere $S^{n-1}\subset\bbr^n$ with
its standard permutation of coordinates action
(see Lemma \ref{deltasphere}). This permutation action
has fundamental domain a simplex. In fact the
quotient $S^{n-1}/K$ is contractible for any $\{1\}\neq
K\subset\mathfrak{S}_n$.  This is readily verified if we think of the
permutation switching the $i$-th and $j$-th coordinate of the sphere as the
reflection with respect to the hyperplane $x_i=x_j$ in $\bbr^n$.
Next and given a simplicial structure
on $X$, there is a simplicial structure on $X^n$ with vertices $\varpi$
given as $n$-tuples
of vertices of $X$. The fat diagonal is sub-simplicial with
vertices having at least two equal entries. Moreover the action of $\mathfrak S_n$
is simplicial and permutes the entries of any $\varpi$ (see \cite{denis}). This simplicial
action on $\Delta_{fat}$ satisfies the hypothesis of proposition 7.11
in \cite{dwyer} and lemma \ref{contractible} follows.
\hfill$\Box$\vskip 5pt

%**************************************************************************

\section{Some Equivariant Deformations}\label{arone}

We clarify and complement some claims we made in the proofs above.

\ble\label{equivariantsplit} Let $(S^1)^{(n)}\cong S^n$ be the
$n$-fold smash product, with $\mathfrak{S}_n$ acting by permuting
factors, and let $\Sigma S^{n-1}$ be the unreduced suspension with
$\mathfrak{S}_n$ acting on the suspension coordinate trivially and
on $S^{n-1}\subset\bbr^n$ by permuting
coordinates fixing the point $p := ({1\over\sqrt{n}},
\ldots, {1\over\sqrt{n}})$. Then both spheres are
$\mathfrak{S}_n$-equivariantly equivalent. \ele

\noindent\proof
Represent $S^1$ as $I=[0,1]$ with $0\sim 1$, so that $(S^1)^{(n)}$ is a
solid cube with all sides identified to a point. We can think of
this cube as inscribed in the unit sphere which via radial stretch
is $\mathfrak{S}_n$-equivariantly identified with the unit disk $D$ in
$\bbr^n$ and the boundary of the cube being identified to this unit sphere.
The action of $\mathfrak S_n$ on $D$ is by reflections with respect
to the hyperplanes $x_i=x_j$. On the other hand, think of the
unreduced suspension $\Sigma S^{n-1}$
as $I\times S^{n-1}/\sim$ with identifications at $t=0,1$.
Let $\mathfrak{S}_n$ act on $I$
trivially and on $S^{n-1}\subset\bbr^n$ by permuting
coordinates.  Then the map $\Phi :
I\times S^{n-1}\lrar D/S$, $(t,x)\mapsto tx$ factorizes through $\Sigma S^{n-1}$
and it is a
$\mathfrak{S}_n$-equivariant equivalence.
\hfill$\Box$\vskip 5pt

\bco\label{suspendaction} Let $Y$ be a pointed left
$\mathfrak{S}_n$-space, and suppose $\mathfrak{S}_n$ acts on
$S^n=(S^1)^{\wedge n}$ by permuting factors. Then
$S^n\wedge_{\mathfrak{S}_n}Y\simeq \Sigma
(S^{n-1}\wedge_{\mathfrak{S}_n}Y)$ and $\mathfrak{S}_n$ acts on
$S^{n-1}\subset\bbr^n$ by
permuting coordinates. In particular we recover the (weaker) equivalence
$$\bsp{n}\Sigma X \simeq \Sigma\sj{n}(X)$$
\eco

\noindent\proof  According to lemma \ref{equivariantsplit},
$S^n\wedge_{\mathfrak{S}_n}Y$ can be replaced by $(S^1\wedge
S^{n-1})\wedge_{\mathfrak{S}_n}Y$ where the symmetric group action
on $S^1$ is trivial and on $S^{n-1}$ is via reflections as above.
The first claim follows. The second claim is obtained through
the series of identifications
\begin{eqnarray*}
\bsp{n}\Sigma X  = (\Sigma X)^{(n)} \simeq S^n
\wedge_{\mathfrak{S}_n} X^{(n)}\simeq S^1\wedge
(S^{n-1}\wedge_{\mathfrak{S}_n}X^{(n)})\simeq \Sigma\sj{n}X
\end{eqnarray*}
$\mathfrak{S}_n$ acting on the smash product by acting
diagonally with action on $S^{n-1}\subset\bbr^n$ as described.
\hfill$\Box$\vskip 5pt

Next we have been using interchangeably the unit sphere in $\bbr^n$ and
the sphere $\Delta/\partial\Delta$ with $\Delta=\Delta_{n-1}$ as in
(\ref{simplex}). This is possible because of the following elementary
observation.

\ble\label{deltasphere}
The sphere $\Delta/\partial\Delta$
is $\mathfrak S_n$-equivariantly homeomorphic to
the unit sphere $S^{n-1}\subset\bbr^n$.
\ele

\noindent\proof Let $H^+\subset S^{n-1}$ be the ``positive hemisphere"
$\{(x_1,\ldots, x_n)\in\bbr^n\ |\ \sum x_i^2=1, x_i\geq 0\ \forall i\}$.
This subspace is invariant under the permutation action of $\mathfrak S_n$.
The point is to see that $H^+$ is a $\mathfrak S_n$ deformation retract
of the unit sphere
$S^{n-1}$, with $S^{n-1}\backslash\{-p\}$ being mapped onto the interior of $H^+$,
 here $p=(1/\sqrt{n},\cdots, 1/\sqrt{n})$. An inverse deformation
 retraction can be given for example by bringing closer to the point $-p$
the ``end points" of the hemisphere- those with a single non trivial entry-
along great circles running through those points, and going from
$p$ and $-p$. Such an inverse retraction extends to $H^+$ in the
obvious way (we get at each time a portion of the sphere containing $p$ and
bounded by
geodesics between the images of these endpoints).
On the other hand the simplex $\Delta$ is  $\mathfrak S_n$-equivariantly
homeomorphic to $H^+$ via the map
$(t_1,\ldots, t_n)\longmapsto \left({t_1},\ldots, t_n\right)/ \sqrt{\sum t_i^2}$.
Composing these equivalences we obtain an equivariant homeomorphism $\Delta \lrar S^{n-1}$, sending the interior of $\Delta$ to $S^{n-1}\backslash\{-p\}$ and
sending the boundary to $-p$. This proves the claim.
\hfill$\Box$\vskip 5pt

\bre\label{vogt} (About the splitting theorem of Vogt).
In \cite{vogt}, it is shown that there is
$\mathfrak{S}_n$-splitting, natural in $X$ up to homotopy and under
$\Sigma^{n-1}D_nX$,
$$\Sigma^{n-1}X^n\simeq \Sigma^{n-1}D_nX\vee X^{*n}$$
where $D_nX\subset X^n$ is the \textit{fat wedge} consisting of all
tuples with basepoint in at least one coordinate.
More precisely, both sides of the equation
are homotopy equivalent through $\mathfrak{S}_n$-equivariant maps.
The case $n=2$ of this theorem is a special
case of the very classical stable splitting of
$\Sigma (X\times Y)$ as the bouquet
$\Sigma X\vee \Sigma Y\vee \Sigma (X\wedge Y)$.
To get the $\bbz_2$-equivariance of the splitting when $X=Y$, the
 involution must act on both $X^2$ and $X\vee X$ by permuting factors, but
also act on the suspension coordinate via $t\mapsto 1-t$. To see why this is so,
we recall that the equivalence $\Sigma (X\times Y)\rightarrow \Sigma X\vee\Sigma Y\vee
\Sigma (X\wedge Y)$ is the wedge of two maps: (i) the projection $\Sigma X\times Y
\rightarrow\Sigma X\wedge Y$, and (ii) the composite of
the pinch map (on the suspension) with the projections
$\pi_1, \pi_2$ of $X\times Y$ onto the first and second factors as in
$$\Sigma (X\times Y)\lrar\Sigma (X\times Y)\vee\Sigma (X\times Y)
\fract{\Sigma (\pi_1\vee\pi_2)}{\ra 4}\Sigma X\vee\Sigma Y$$
The point is that in order for this map to commute with the switch map interchanging
$X$ and $Y$, one needs to switch the suspension
coordinate from $t$ to $1-t$ as well.
\ere

%***************************************************************

\section{Connectivity}\label{connectivity}

The aim of this section is to prove Theorem \ref{cor1}. All spaces
$X$ are assumed again to be connected and based simplicial
complexes. We first establish a few basic facts about relative
Leray-Serre spectral sequences for ``singular'' maps; that is
maps $E\rightarrow X$ which might fail to be fibrations
over some closed subspace $A\subset X$. This situation occurs when for example we
collapse out some fibers from the total space of a fibration, or if we have
the quotient map of a group action which is free away from some fixed set $A$.
Such a situation is ubiquitous in topology.

\bde\label{hypothesis} We say $(F,F_0)\lrar (E,W)\fract{f}{\lrar} X$
is a fibration \textit{relative to $A\subset X$} if\ : \\
(i) $A$ is closed and $f$ is a fibration away from $A$; that is
$$(F,F_0)\lrar
(E-f^{-1}(A), W-f^{-1}(A)\cap W)\lrar X-A$$ is a fibration
(typically in the examples below a bundle). This means that
$(F,E-f^{-1}(A),X-A)$ is a fibration of which $(F_0, W-f^{-1}(A)\cap W,X-A)$ is
a subfibration.\\
(ii) There is an open neighborhood $A\subset V$
such that $A$ is a deformation retract of $V$ and the deformation
retraction lifts upstairs; that is the pair ($f^{-1}(V),
f^{-1}(V)\cap W)$ deformation retracts onto ($f^{-1}(A),
f^{-1}(A)\cap W$).
\ede

\ble\label{key1} Let $(F,F_0)\lrar (E,W)\lrar X$ be a fibration
relative to $A\subset X$, with $X$ simply connected.
Then there is a spectral sequence with
$$E^2 = H_*(X,A; H_*(F,F_0))\Longrightarrow
H_*(E,W\cup f^{-1}(A))$$ \ele

\noindent\proof   We have a
fibration $(F,F_0)\lrar (E',W')\lrar X-A$ where $E'=E-f^{-1}(A)$ and
$W'=W-f^{-1}(A)\cap W$. By the hypothesis we can choose a
neighborhood $V$ of $A$ so that both $V$ and its closure $\overline{V}$ deformation
retract onto $A$ while $f^{-1}(V)$ and $f^{-1}(\overline{V})$ deformation
retract onto $f^{-1}(A)$. Now $\overline{V}-A$ is closed in $X-A$ and
according to \cite{morin} there is a spectral sequence
$$E^2 = H_*(X-A,\overline{V}-A;  H(F,F_0)) \ \Longrightarrow\
H_*(E-f^{-1}(A),f^{-1}(\overline{V}-A)\cup W')$$ We can replace
$H_*(X-A,\overline{V}-A)\cong H_*(X,\overline{V})$ by excision and then
$H_*(X,\overline{V})\cong H_*(X,A)$ since $\overline{V}$ deformation retracts
onto $A$. On the other hand
\begin{eqnarray*}
H_*(E-f^{-1}(A),f^{-1}(\overline{V}-A)\cup W')&\cong&
H_*(E-f^{-1}(A),f^{-1}(\overline{V}-A)\cup W-f^{-1}(A)\cap W)\\
&\cong& H_*(E-f^{-1}(A), f^{-1}(\overline{V})\cup W -f^{-1}(A))\\
& \cong& H_*(E, f^{-1}(\overline{V})\cup W)\\
&\cong &H_*(E,f^{-1}(A)\cup W)
\end{eqnarray*}
The one before last isomorphism is obtained by excision again ($f^{-1}(A)$
being closed), and the last isomorphism follows upon comparing the long
exact sequences for $(f^{-1}(\overline{V})\cup W, W)$ and
$(f^{-1}(A)\cup W,W)$, and using the fact that ($f^{-1}(\overline{V}),
f^{-1}(\overline{V})\cap W)$ deformation retracts onto ($f^{-1}(A),
f^{-1}(A)\cap W$).
\hfill$\Box$\vskip 5pt

\bco\label{key2} Let $F\lrar E\fract{f}{\lrar} X$ be a fibration
relative to $A$ which admits a section $s : X\rightarrow E$. We view
$X$ as a subspace of $E$ through the section and assume it is simply
connected. Then there is a spectral sequence converging to
$H_*(E, f^{-1}(A)\cup X)$ with $E^2$-term $H_*(X,A; \tilde H_*(F))$.
\eco

The proof is a direct application of lemma \ref{key1} applied to
$(F,*)\lrar (E,X)\lrar X$ which is a fibration relative to $A$. We
are now in the position to prove Theorem \ref{cor1} of the
introduction.

\bth\label{connectivitybound} If $X$ is $r$-connected, $r\geq 1$,
then ${\mathcal B}_n(X)$ is $2n+r-2$-connected.
\end{theorem}

\noindent\proof
Write ${\mathcal B}_n(X)\simeq
S^{n-1}\wedge_{\mathfrak{S}_n}(X^{(n)}/\Delta_{fat})$ as in the proof of
Theorem \ref{equivalence}.  Notice that the action of
$\mathfrak{S}_n$ on $X^{(n)}/\Delta_{fat}$ is free everywhere but at
the basepoint $*$ to which $\Delta_{fat}$ is collapsed out. This is
clear since fixed points of $\mathfrak{S}_n$ acting on $X^n$ consist
precisely of tuples with two or more equal entries.

Let $E_n=S^{n-1}\times_{\mathfrak{S}_n}(X^{(n)}/\Delta_{fat})$ and
set $X_n:=*\times_{\mathfrak{S}_n}(X^{(n)}/\Delta_{fat}) =
\left(X^{(n)}/\Delta_{fat}\right)_{\mathfrak{S}_n}$ where $* =
(1/\sqrt{n},\ldots, 1/\sqrt{n})$ is the fixed point under the
$\mathfrak{S}_n$-action. Let $x_0$ denote both the unique fixed
point of $X^{(n)}/\Delta_{fat}$ under the $\mathfrak{S}_n$-action
and its image in the quotient $X_n$. The projection $f : E_n\lrar
X_n$ has as preimage $S^{n-1}$ over every point of $X_n-\{x_0\}$ and over
$x_0$ the fiber is $S^{n-1}/\mathfrak{S}_n$ which is contractible as
we pointed out earlier. We claim that $f$ is a bundle relative to $x_0$.
This is straightforward. Let $V$ be a neighborhood of $x_0\in X_n$
that deformation retracts onto it with deformation retraction $G$.
Write $p_n: X^{(n)}/\Delta_{fat}\lrar X_n$ the projection. Then
$p_n^{-1}(V-\{x_0\})\lrar V-\{x_0\}$ is a regular covering with a
single ramification point at $p_n^{-1}(x_0)=x_0$. The deformation
$G$ of $V$ onto $x_0$ lifts to a $\mathfrak{S}_n$-equivariant
deformation of $p_n^{-1}(V)$ onto $x_0$ and hence lifts to a
deformation of $f^{-1}(V) = S^{n-1}\times_{\mathfrak{S}_n}
p_n^{-1}(V)$ onto $f^{-1}(x_0)$. Our hypothesis in definition
\ref{hypothesis} for $f$ as above, $(F,F_0)=(S^{n-1},*)$ and $A=x_0$ are
verified.

The space $X_n$ has been studied by Nakaoka where he has shown
(\cite{nakaoka}, Proposition 4.3) that if $X$ is $r$-connected, then
$X_n=(X^{(n)}/\Delta_{fat})_{\mathfrak{S}_n}$ is $r+n-1$ connected and
this only holds for $r\geq 1$ hence our hypothesis\footnote{It seems
very likely that we can weaken the condition $r\geq 1$ to $r\geq 0$ if
$X$ is a manifold.}. In
particular $X_n$ is always simply connected if $X$ is.

Notice that $E_n\lrar X_n$ has a section and so according to
corollary \ref{key2} we get a spectral sequence
\begin{equation}\label{e2}
E^2 = \tilde H_*(X_n, \tilde H_*(S^{n-1}))\ \Longrightarrow\
H_*(E,f^{-1}(x_0)\cup X_n)\cong H_*({\mathcal B}_n(X))
\end{equation}
This $E^2$-term is necessarily trivial in total degree $r+n-1+n-1 =
2n+r-2$ using Nakaoka's aforementioned connectivity result so that
$\tilde H_q({\mathcal B}_n(X))=0$ for $q\leq 2n+r-2$. Since ${\mathcal B}_n(X)$ is
simply connected (Lemma \ref{simplyconnected}), the claim follows.
\hfill$\Box$\vskip 5pt

\bre Since $\Sigma {\mathcal B}_n(X)\simeq\bsp{n}(\Sigma X)$ (Proposition
\ref{key}) the above calculation shows that the connectivity of
$\bsp{n}(\Sigma X)$ is $2n+r-1$ if $X$ is $r$-connected. In fact we
can deduce from here that $\bsp{n}(X)$ is $2n+r-2$ connected
if $r\geq 1$ is the connectivity of $X$ \cite{braids}. \ere

%***************************************************************

\section{The Space of Chords}\label{chords}

Given a topological space $X$, then ${\mathcal B}_2(X)$ is referred to as the
space of chords on $X$ and can be viewed as the space of all
segments $pq = qp$ between any two point $p,q\in X$. The length of
$pq$ is assumed to be a continuous function of $p$ and $q$ so that
the length approaches zero if $p$ and $q$ approach a common limit.
Distinct segments do not intersect except at a common endpoint. This
space has been looked at by Clark in 1944 \cite{clark} who wrote
down a homological description in terms of simplicial generators and
relations. The method of Clark consisted in analyzing a
Mayer-Vietoris exact sequence and is hard to work with.
Below we use Proposition \ref{key} and
the identity
\begin{equation}\label{decompose}
\bsp{n}(X_1\vee\cdots\vee X_k)\cong\bigvee_{r_1+\cdots +r_k=n}\bsp{r_1}X_1\wedge
\cdots\wedge\bsp{r_k}X_k
\end{equation}
to analyze some useful cases.

\ble\label{product} For based spaces $X$ and $Y$, ${\mathcal B}_2(X\times Y)$
is stably homotopy equivalent to a wedge of six terms
$${\mathcal B}_2(X)\vee {\mathcal B}_2(Y)\vee {\mathcal B}_2(X\wedge Y )
\vee (X*  Y)\vee (X^{*2}\wedge Y)\vee (Y^{*2}\wedge X)$$ the
equivalence occurring after a single suspension.\ele

\noindent\proof  The claim comes down to decomposing the
symmetric smash
\begin{eqnarray*}
\Sigma {\mathcal B}_2(X\times Y) &=&
\bsp{2}(\Sigma (X\times Y))\\
&\simeq&\bsp{2}(\Sigma X\vee\Sigma Y\vee (\Sigma X\wedge Y))\\
&\simeq&\bsp{2}(\Sigma X)\vee\bsp{2}(\Sigma Y)\vee\bsp{2}(\Sigma X\wedge Y)
\vee \Sigma^2(X\wedge Y)\vee \Sigma^2(X\wedge X\wedge Y)\vee
\Sigma^2(Y\wedge X\wedge Y)
\end{eqnarray*}
and this is our claim after desuspending once.
\hfill$\Box$\vskip 5pt

Let $C_g$ be a closed topological surface of genus $g$ and write as before
$X^{\vee n}$ the one-point union of $n$-copies of $X$. Then

\ble\label{closedsur} $\Sigma {\mathcal B}_2(C_g)\simeq (S^4)^{\vee
(2g^2+g)}\vee (S^5)^{\vee 2g}\vee\Sigma^4\bbr P^2$.  \ele

\noindent\proof
The surface $C_g$ is obtained by attaching a cell of dimension two
to a bouquet of $2g$
circles $\bigvee^{2g}S^1$. The suspension of the attaching map is trivial
so that we have a standard splitting
$\Sigma C_g\simeq S^3\vee\bigvee^{2g}S^2$.
We can then replace the symmetric smash
$\bsp{2}(\Sigma C_g)$ by
$\bsp{2}\left(S^3\vee\bigvee^{2g}S^2\right)$ and apply (\ref{decompose}).
\hfill$\Box$\vskip 5pt

\bex As it should be, both lemmas \ref{product} and
\ref{closedsur} agree for the case of the torus $S^1\times S^1$ and
one has ${\mathcal B}_2(S^1\times S^1)\simeq_s S^3\vee S^3\vee S^3\vee
S^4\vee S^4\vee\Sigma^3\bbr P^2$.  \eex

%*******************************************************************

\section{Barycenter Spaces of Spheres}\label{spheres}

The main objective of this section is to show that ${\mathcal B}_n(S^k)\simeq
\sj{n}(S^k)$ is a $(k+1)$-fold suspension on an explicit space.  The special
cases of $n=2$ and of the circle were discussed in lemma \ref{jmea1}
and in corollary \ref{circle} respectively. As a start we recall some useful
constructions from \cite{su}.

There is a homeomorphism $S^{n}* S^{m}= \partial (D^{n+1}\times
D^{m+1})=S^{n+m+1}$ (see examples \ref{examples}, (1)).
In fact consider
$S=S^{n+m+1}$ the unit sphere in $\bbr^{n+m+2}$ and choose an
orthogonal decomposition $\bbr^{n+1}\oplus\bbr^{m+1}$ so that $S^n=
S \cap\bbr^{n+1}$, $S^m=S\cap\bbr^{m+1}$. Then by joining segments
between these two spheres we obtain $S^n* S^m$ and radial
projection onto $S$ yields the homeomorphism
\begin{equation}\label{specialhom}
\Phi : S^n*  S^m\lrar S^{n+m+1}\ ,\ [x, y, t]\mapsto
\sqrt{1-t}x+\sqrt{t}y
\end{equation}

Slightly more generally let $D^l$ be the unit disk in $\bbr^l$ and
$S^{l-1}$ its boundary. We write $(D,S)$ for the corresponding pair
and we usually omit to specify $l$ when it is clear from the
context. Suppose $W_1\oplus W_2=\bbr^l$ is an orthogonal
decomposition and set $D_i=W_i\cap D$, $S_i=W_i\cap S$. Then if $*$
is the join product as before, we have a homeomorphism of pairs
\begin{equation}\label{firsthomeo}
\Phi : (D_1*S_2, S_1*S_2)\fract{\cong}{\lrar} (D,S)
\end{equation}

We set $l=nk$ and $D=D^{nk}\subset (\bbr^k)^n$ the unit disk. There is an
action of $\mathfrak{S}_n$ on $(\bbr^k)^n$  by permuting vectors.
This action is via orthogonal transformations and descends to an
action on the pair $(D,S)$. Consider the $\mathfrak{S}_n$-invariant
subspace
\begin{equation}\label{thespaceh}
W_1=\{{\bf v}=(v_1,\ldots, v_n)\in (\bbr^k)^n\ |\ v_1=\ldots = v_n\}
\end{equation}
Its orthogonal complement $W_2=W_1^{\perp}$ is
$\mathfrak{S}_n$-invariant, and $\mathfrak{S}_n$ acts on the joins
$D_1*S_2$ and $S_1*S_2$ by acting diagonally on factors. Since the
action on $(D_1,S_1)$ is trivial, and since the homeomorphism
(\ref{firsthomeo}) is $\mathfrak{S}_n$-equivariant, one obtains
\begin{equation}\label{same}
(D^{nk}/\mathfrak{S}_n, S^{nk-1}/\mathfrak{S}_n )\cong
(D_1*(S_2/\mathfrak{S}_n), S_1*(S_2/\mathfrak{S}_n))
\end{equation}
This leads to the interesting consequence \cite{su}.

\bpr\label{su} $\bsp{n}S^k = S^k*(S_2/\mathfrak{S}_n) =
\Sigma^{k+1}(S_2/\mathfrak{S}_n)$.
\epr

\noindent\proof
As in \cite{su}, p. 370,
the pair $(D,S)$ in $(\bbr^k)^n$ is $\mathfrak{S}_n$-equivariantly
equivalent to the pair
$$((D^k)^n, \partial (D^k)^n)\ \ \  \ \ \hbox{where}\ \ \ \partial
(D^k)^n:= \bigcup (D^k)^i\times S^{k-1}\times (D^k)^{i-1}$$
The quotient
$(D^k)^n/\partial (D^k)^n$ is $(S^k)^{(n)}$ and the permutation
action of $\mathfrak{S}_n$ translates into a permutation of the
smash factors so that the quotient is $\bsp{n}S^k$. But according to
(\ref{same}) this quotient must coincide with
$$[D_1*(S_2/\mathfrak{S}_n)]/[S_1*(S_2/\mathfrak{S}_n)]
\simeq (D_1/S_1)*(S_2/\mathfrak{S}_n) = S^k*(S_2/\mathfrak{S}_n)$$ where
again $D_1$ is a disc of dimension $k=\dim W_1$.
\hfill$\Box$\vskip 5pt

We have just shown that $\Sigma \sj{n}(S^{k-1}) = \bsp{n}(S^k)$ is
a $k+1$-fold suspension.
The main statement next is that this identification desuspends.
With $S_2$ as above the unit $({k(n-1)-1})$-sphere in $W_2=W_1^\perp$, we have

\bth\label{keyresult} There is a homotopy equivalence
$\sj{n}(S^{k-1})\simeq\Sigma^{k}\left(S_2/{\mathfrak{S}_n}\right)$.
\end{theorem}

\noindent\proof We first adopt
some notation. Define the smash product quotient
\begin{equation}\label{redjoin}
L_n(X):= S^{n-1}\wedge_{\mathfrak{S}_n}X^{(n)}
\end{equation}
where $S^{n-1}$ is $\Delta/\partial\Delta$ with
the appropriate $\mathfrak S_n$-action.
Of course $L_n(X)\simeq\sj{n}(X)$ according to Theorem \ref{equivalence}.
Write $Y_{\mathfrak{S}_n}$ for the quotient of $Y$ by a
$\mathfrak{S}_n$ action. If $A$ and $B$ are compact manifolds with boundary
$\partial$, define as well
$$A\wedge^{\partial} B:= A\times B/\partial A\times
B\cup A\times\partial B = (A/\partial A)\wedge (B/\partial B)$$ For
example $I\wedge^{\partial} I = S^2$. If moreover $G$ is a group
acting on the right of $A$ and on the left of $B$ preserving
boundaries, then we write $A\wedge^{\partial}_GB$ for the quotient
of $A\wedge^{\partial}B$ under the diagonal action.
For instance $\Delta:=\Delta_{n-1}\subset I^n$
and $\Delta\times (D^{k-1})^n$ is a $\mathfrak{S}_n$-invariant
subspace of $I^n\times (D^{k-1})^n$ (under the diagonal action) so
that we can write
$$L_n(S^{k-1})=
S^{n-1}\wedge_{\mathfrak{S}_n}(S^{k-1})^{(n)}
= \Delta\wedge^{\partial}_{\mathfrak{S}_n}(D^{k-1})^n
$$
We will transform this quotient in steps. We will use freely the
identifications (equivariant with respect to the permutation action)
$[0,1]^n\cong [-1,1]^n\cong D^n$. Also $(D^k)^n$ and $D^{kn}$ are
$\mathfrak S_n$-equivariant subspaces of $(\bbr^k)^n$. Let now
\begin{equation}\label{vn}
V_n = \{(t_1,\ldots, t_n)\in [-1,1]^n\ |\ \sum t_i = 0\}
\end{equation}
and consider the sequence of $\mathfrak S_n$-equivariant homeomorphisms
\begin{eqnarray}
\Delta\times (D^{k-1})^n&\cong&
\Delta\times ([-1,1]^{k-1})^n\nonumber\\
&\cong&V_n\times ([-1,1]^{k-1})^n\nonumber\\
&\cong&\{(v_1,\ldots, v_n)\in ([-1,1]^k)^n\ |\ \sum v^1_i= 0\}\label{id2}\\
&\cong&\{(v_1,\ldots, v_n)\in D^{kn}\ |\ \sum v^1_i= 0\} := B_n \label{bn}
\end{eqnarray}
where if $v\in [-1,1]^k$, then $v^1\in [-1,1]$ is its first
coordinate. The identification (\ref{id2})
is obtained by distributing coordinates in the obvious way. To see how to
get (\ref{bn}) we need check that the equivalences $([-1,1]^k)^n\cong
(D^k)^n\cong D^{kn}$
preserve the subspace $\sum_{i=1}^n v^1_i=0$. This is clear for the first
 equivalence. To see the second, we could observe that $(D^k)^n$ is the unit sphere
 in $(\bbr^{k})^n$ but for a different norm. For $v=(v_1,\ldots, v_n)\in (\bbr^{k})^n$ define (as in \cite{su}, p.370)
$$||v||'= max_i||v_i||$$
Then $(D^k)^n$ is the subspace of all such $v$ with $||v||'=1$.
The map $v\mapsto v'={||v||\over ||v||'}v$ maps $D^{kn}$ homeomorphically and
equivariantly onto $(D^k)^n$. Obviously if $\sum v^1_i=0$, then $\sum v'^1_i=0$
as well. This establishes the last equivalence (\ref{bn}).

Note that the space $B_n$ in (\ref{bn}) is the unit sphere in the hyperplane
$\{(v_1,\ldots, v_n)\ |\ \sum v_1^i= 0\}$.
It is identified equivariantly with $\Delta\times
(D^{k-1})^n$. Since this identification preserves boundaries we see
that
$$L_n(S^{k-1}) = \Delta\wedge^{\partial}_{\mathfrak{S}_n}(D^{k-1})^n
\cong (B_n/\partial B_n)_{\mathfrak{S}_n}$$ The problem then boils down
to identifying the $\mathfrak{S}_n$-equivariant homotopy type
of the pair $(B_n,\partial B_n)$. This is where the methods of
proposition \ref{su} come handy.

As in (\ref{thespaceh}), set
$W_1=\{({\bf v}\in (\bbr^k)^n\ |\
v_1=\ldots = v_n\}$ and the various disks
$$D=D^{nk}\ \ ,\ \ D_1=D\cap W_1\ \ ,\ \
D'_1=B_n\cap W_1$$ and $D_1^{\perp}$ the unit disk in $W_1^{\perp}$.
The associated boundary spheres are denoted by $S, S_1, S'_1$ and
$S_1^{\perp}$ respectively. Since $W_1^{\perp} = \left\{{\bf v} =
(v_1,\ldots, v_n)\in (\bbr^k)^n\ |\ \sum v_i = 0\right\}$, the condition
$\sum v_i^1=0$ says that $S_1^{\perp}\subset D_1^{\perp}\subset
B_n$. Consider the homeomorphism
$\Phi : (D_1*S_1^{\perp},S_1*S_1^{\perp})\rightarrow
(D, S)$ constructed by means of the map (\ref{specialhom}).
 Its restrition $\Phi_{|}$ makes the following diagram commute
$$\xymatrix{
(D'_1*S_1^{\perp}, S'_1*S_1^{\perp})\ar[d]^\subset\ar[r]^{\ \ \ \
 \Phi_{|}}&(B_n,\partial B_n)\ar[d]^\subset\\
(D_1*S_1^{\perp},S_1*S_1^{\perp})\ar[r]^{\ \ \ \ \Phi}& (D, S)
}
$$
and $\Phi_{|}$ is a homeomorphism of pairs as well (this is easy
to see since $\Phi$ restricted to $D'_1*S_1^{\perp}$
surjects onto $B_n$).
We are at last ready to conclude. We have a series of equivalences
$$B_n/\partial B_n \ \cong\
D'_1*S_1^{\perp}/S'_1*S_1^{\perp}\ \simeq\ (D'_1/S'_1)*S_1^{\perp} =
S^{k-1}*S_1^{\perp}\ \ \ ,\ \ \ \dim D'_1 = k-1$$ The action of
$\mathfrak{S}_n$ is trivial on $S^{k-1}$ since it is trivial on
$D_1$ and thus on $D'_1$, so by passing to quotients
$$L_n(S^{k-1}) = (B_n/\partial B_n)_{\mathfrak{S}_n}\
\simeq\ S^{k-1}*(S_1^{\perp}/\mathfrak{S}_n )\ \simeq\
\Sigma^k(S_1^{\perp}/\mathfrak{S}_n )$$ and the proof is complete.
\hfill$\Box$\vskip 5pt

%*****************************************************************

\section{The case of Manifolds}\label{manifolds}

We shed different light on calculations of Bahri-Coron and Lannes
\cite{bahri} on the homology of barycenters of manifolds. First we
have the following description of the top homology group for
symmetric products of general closed manifolds.

\bpr\label{topcell}\cite{denis}
Suppose $M$ is a closed manifold of dimension $d\geq 2$.
If $M$ is orientable, then
$$
H_{nd}(\sp{n}M;\bbz ) =
\begin{cases}\bbz, &d\ \hbox{even}\\
0,&d\ \hbox{odd}
\end{cases}
$$
For general closed manifolds $M$,
$H_{nd}(\sp{n}M;\bbz_2 ) = \bbz_2$.
\epr

\bco\label{topdim} If $M$ is closed orientable of dimension $d\geq 1$, then in
top dimension $H_{n(d+1)-1}({\mathcal B}_n(M)) = \bbz$ if $d$ is odd and is
$0$ if $d$ is even. \eco

\noindent\proof  By Theorem \ref{main2} we need determine
$H_{n(d+1)}(\bsp{n}(\Sigma M))$. Since $\Sigma M$ is a CW complex
with top integral homology group $H_{d+1}=\bbz$; $\Sigma M$ has the
homology of a wedge $S^{d+1}\vee Y$ where $Y$ is of dimension $d$. By a
well-known result of Dold, the homology of symmetric products of
simplicial complexes, and hence of their reduced symmetric products,
only depends on the homology of the underlying space. This means
that
$$H_*(\bsp{n}(\Sigma M))\cong
H_*(\bsp{n}(S^{d+1}\vee Y)) \cong \bigoplus_{r+s=n}
H_*(\bsp{r}S^{d+1}\wedge \bsp{s}Y)$$ according to the formula (\ref{spwedge}).
Since the homological dimension
of $Y$ is $d$, the term
$H_*(\bsp{r}S^{d+1}\wedge \bsp{s}Y)$ is trivial in degrees larger than
$r(d+1) + sd$. It follows that in top dimension
$$H_{n(d+1)}(\bsp{n}(\Sigma M))\cong H_{n(d+1)}(\bsp{n}S^{d+1}) \cong
H_{n(d+1)}(\sp{n}S^{d+1})$$
The claim now follows from Proposition \ref{topcell} applied to
$M=S^{d+1}$.
\hfill$\Box$\vskip 5pt

\bex ${\mathcal B}_2(S^k)\simeq\Sigma^{k+1}\bbr P^{k}$ by (\ref{jmea}) and
hence $H_{2(k+1)-1}({\mathcal B}_2(S^k))=H_{k}(\bbr P^k)$ and this is indeed
$\bbz$ or $0$ according to whether $k$ is odd or even. Similarly and
by lemma \ref{closedsur}, the top class in dimension 5 of the space
of chords of a closed Riemann surface is trivial as expected
$$H_5({\mathcal B}_2(C_g))\cong H_6((S^4)^{\vee
(2g^2+g)}\vee (S^5)^{\vee 2g}\vee\Sigma^4\bbr P^2) =
H_6(\Sigma^4\bbr P^2) = 0$$\eex

In \cite{bahri} a main point of consideration were tranfer morphisms
$$\Phi : H_{nd+n-1}({\mathcal B}_nM,{\mathcal B}_{n-1}M)\lrar
H_{(n-1)d+n-2}({\mathcal B}_{n-1}M,{\mathcal B}_{n-2}M)$$ sending orientation class to
orientation class. In this case $M$ was allowed to have boundary and
coefficients were in $\bbz_2$. The map $\Phi$ was a mix of a
transfer map, a cap product and a boundary morphism. We indicate
below a streamlined construction of a transfer map which appeals as
before to the identification of $\Sigma {\mathcal B}_n(M)$ with
$\bsp{n}(\Sigma M)$. In our case and for closed oriented $M$ we seek
to construct for each positive $n$ a map
\begin{equation}\label{trans}\Theta : H_{n(d+1)}(\bsp{n}(\Sigma M))\lrar
H_{(n-1)(d+1)}(\bsp{n-1}(\Sigma M)) \end{equation}
with the right
homological properties. The construction is due to L. Smith and has
refinements in \cite{dold}. Consider the degree $n$ covering
\begin{eqnarray}\label{smithtransfer}
\Sigma M\times \sp{n-1}(\Sigma M)&\fract{p}{\lrar}&\sp{n}(\Sigma M)\\
(x, y_1+\cdots + y_{n-1})&\longmapsto& x+y_1+\ldots + y_{n-1}\nonumber
\end{eqnarray}
The transfer map in homology for this covering can be achieved at
the level of spaces as follows. Set $Y=\Sigma M$ to ease notation.
To each element $y_1+\ldots + y_n\in \sp{n}(Y)$, we associate the
unordered tuple
\begin{equation}\label{derivative}
\sum (y_i, y_1+\ldots + y_{i-1}+y_{i+1} + \ldots+ y_n)\ \ \in\ \
\sp{n} (Y\times \sp{n-1}Y)
\end{equation}
The induced map in homology is a map
$\tau : H_*(\sp{n}Y)\lrar H_*(\sp{n}(Y\times\sp{n-1}Y))$.

Next we invoke a homological splitting of Steenrod which asserts
that for any based connected space $X$, $H_*(X,\bbz )$ is a \textit{direct
summand} of $H_*(\sp{n}X)$ (see (\ref{steensplit})). We then get
a projection which we write $st: H_*(\sp{n}X)\lrar H_*(X)$.
We can now define our transfer map $\Theta$ as the composite
\begin{eqnarray*}
\Theta :
H_*(\sp{n}Y)\fract{\tau}{\lrar} H_*(\sp{n}(Y\times \sp{n-1}Y))
&\fract{st}{\lrar}& H_*(Y\times\sp{n-1}Y)\\
&\lrar&H_*(S^{d+1}\times\sp{n-1}Y)\\
&\fract{int}{\lrar}& H_*(S^{d+1}\wedge\sp{n-1}Y)\cong
H_{*-d-1}(\sp{n-1}Y)
\end{eqnarray*}
where $Y=\Sigma M\lrar S^{d+1}$ is the suspension of the map
$M\lrar S^d$ which collapses the complement of an open disk in $M$.
Write $v_1\in
H_{d+1}(Y;\bbz_2 )$ this top class, and let
$v_n$ be the top class in $H_{n(d+1)}(\sp{n}Y;\bbz_2 )\cong H_{n(d+1)}(\sp{n}S^{d+1})$.
Write $tr := st\circ \tau$ the composite of the first two maps making
up $\Theta$. Then

\ble\label{transferclass} $tr (v_n) = v_1\tensor v_{n-1}$
and $\Theta_*(v_n)=v_{n-1}$. \ele

\noindent\proof
The map $tr: H_*(\sp{n}Y)\lrar H_*(Y\times\sp{n-1}Y)$ is a \textit{transfer
  map} for the degree $n$ branched covering (\ref{smithtransfer}), that is the
  composite $p_*\circ tr$ is multiplication by $n$ in homology (see \cite{dold}).
On the other hand the
projection $Y\times\sp{n-1}(Y)\lrar\sp{n}(Y)$ is an $n$-fold branched covering
and hence on top dimension it is multiplication by $n$. It follows that
$$p_*(tr  (v_n))
= nv_n = p_*(v_1\tensor v_{n-1})$$
from which we deduce that $tr
(v_n) = v_1\tensor v_{n-1}$ since  $v_1\tensor
v_{n-1}$ is the only generator in that dimension.
The equality $\theta_*(v_n)=v_{n-1}$ is immediate.
\hfill$\Box$\vskip 5pt

It would be a good exercise to
check whether this transfer agrees or not with the transfer map
constructed in appendix C of \cite{bb, bahri}, but we don't pursue this
further.

%********************************************************************

\section{Appendix : Homology Computations}

This final appendix summarizes for convenience some known results
about the homology of the reduced symmetric products $\bsp{n}X$ and
derives information about the homology of barycenter spaces of spheres
and manifolds.

As is standard $\spy (X)$ will denote
the direct limit under the inclusions $\sp{n}( X) \lrar\sp{n+1}( X)$,
$\sum x_i\mapsto \sum x_i+x_0$ where $x_0\in X$ is a chosen basepoint
($X$ always assumed to be connected).
Since $\spy (X)$, is a connected
abelian topological monoid (associative with unit), it has the homotopy
type of a product of Eilenberg-MacLane spaces. In fact and as a consequence
of a theorem of Dold and Thom
$$\spy (X)\simeq\prod_i K(\tilde H_i(X),i )$$
where $K(G,i)$ is the standard Eilenberg-MacLane space with
$\pi_i(K(G,i))=G$ (for a nice account see \cite{spanier}).
On the other hand and by a computation of
Steenrod (see \cite{dold2}) there is a \textit{splitting}
\begin{eqnarray}\label{steensplit}
H_*(\sp{n}X)&=& H_*(\sp{n-1}X)\oplus H_*(\sp{n}X,\sp{n-1}X)\\
&=& H_*(\sp{n-1}X)\oplus \tilde H_*(\bsp{n}X )\nonumber
\end{eqnarray}
so that $H_*(\bsp{n}X)$ embeds as a direct summand in
$H_*(\sp{n}X)$. Replacing $X$ by a suspension $\Sigma X$ in the expressions
above yield

\bco $H_*({\mathcal B}_n(X))$ embeds as a direct summand in $H_{*+1}(\sp{n}\Sigma X)$
which
in turn embeds in $\bigotimes_{i=0} H_{*+1}(K(\tilde H_i(X),i+1))$.
\eco

We propose to understand the image of this embedding for spheres.
This follows principally from work of Nakaoka and Milgram
\cite{jim}. Since $\spy (S^n)$ is a $K(\bbz, n)$, we can use the
Steenrod splitting above and write
\begin{equation}\label{splitsphere}
\tilde H_*(K(\bbz, n);\bbf )\cong
\bigoplus_{j=1}\tilde H_*(\bsp{j} S^n,\bbf )
\end{equation}
This splitting is valid in cohomology (additively) as well.
The next step is to filter $H_*(K(\bbz, n);\bbf )$ over the positive
integers so that $H_*(\bsp{j}S^n;\bbf )$ corresponds to the classes
of filtration degree precisely $j$. This inductive procedure yields
the following computations which can be deduced from \cite{jim,nakaoka2}.

Recall by work of Serre (Theorem 3 of \cite{serre}) that for $n\geq 2$,
$H^{*}(K(\bbz, n);\bbf_2 )$ (resp. $H^{*}(K(\bbz_2,n);\bbf_2 )$)
is a vector space having as basis elements all iterated Steenrod squares
$Sq^{i_1}\cdots Sq^{i_r}u$ on an $n$ dimensional class $u$ and running
 over sequences of strictly positive integers
$I=\{i_1,\ldots, i_r\}$ which satisfy the
\textit{admissibility} conditions
(i), (ii) and (iii) below (resp. (i) and (ii)):\\
(i) $i_1-i_2\cdots -i_r<n$,\\
(ii) $i_{k}\geq 2i_{k+1}$, $k=1,2,\ldots, r-1$,\\
(iii) $i_r > 1$.

The graded group
$H^*(\sp{k}S^n;\bbf_2)$ is a direct summand of $H^{*}(K(\bbz, n);\bbf_2 )$ as
we previously explained. To describe it,
let $\iota_n\in H^n(\sp{k}S^n)$ be the bottom generator, and set
the \textit{filtration degree}
of $Sq^I(\iota_n )$ to be $2^r$ where $I=\{i_1,\ldots, i_r\}$ as above. Then

\bth\label{jim} For $n\geq 2$.
$H^*(\sp{k}S^n;\bbf_2 )$ is the quotient of the polynomial algebra
$\bbz_2[Sq^I(\iota_n )]$, $I$ running over admissibles, by the ideal
of elements having filtration degree greater than $k$. The subspace
$H^*(\bsp{k}S^n;\bbf_2 )$ corresponds to elements of exact filtration $k$.
\end{theorem}

\noindent\proof
This is a reformulation of the calculation of Nakaoka \cite{nakaoka2, ucci}.
\hfill$\Box$\vskip 5pt

A similar result holds modulo $p$ for odd $p$ with filtration degree
set at $p^r$. Note that
the results in \cite{jim} were stated in homology and the above is the
dual version. A more geometrical viewpoint on symmetric products of
spheres can be found in work of Ucci \cite{ucci}.

For a filtered
$H$-space $X$, $A:=H_*(X;\bbf_p)$ becomes a bigraded algebra over
$\bbf_p$ and so we write $x_{(i,k)}\in A$ for an element $x$ of
\textit{homological} degree $i$ and \textit{filtration} degre $k$. Both the homology
and filtration degrees are additive under the product so that
$x_{(i,k)}\cdot y_{(j,s})$ is a class of degree $i+j$ and filtration
$k+s$.

\bex\label{123} Since $\spy (S^1)\simeq S^1$,
$H^*(\spy (S^1);\bbz ) = E[e_{(1,1)}]$ where $E$ stands for exterior
algebra. On the other hand,
$H^*(\spy (S^2);\bbz )\cong \bbz [b_{(2,1)}]$ and this is in turn consistent with
the diffeomorphism $\spy (S^2)\cong\bbp^\infty$. The case of $\spy (S^3)$
is much more delicate and is given with mod-$2$ coefficients by
$$H^*(\spy (S^3),\bbf_2)\cong \bbf_2[\iota_{(3,1)},
f_{(5,2)}, f_{(9,4)},\ldots, f_{(2^{i}+1,2^i)},\ldots ]$$
where $f_{(2^{2i+1},2^i)}:= Sq^{2^i}\cdots Sq^4Sq^2(\iota )$ with $\iota\in H^3(\spy S^3)$
is the bottom class. A similar formula mod-$p$ is written up in
\cite{c2m2} (see corollary \ref{homcalc} below).
\eex

For $A$ a bigraded algebra as above, we will denote by
$[A]_k$ the submodule of elements with filtration degree precisely $k$.
Let $\sigma$ be the formal suspension operator which raises
homological degree by one. If $E[-]$ stands for an exterior
algebra over the field with $p$ elements $\bbf_p$, then the following
is a consequence of Theorem \ref{main2} and example \ref{123}.

\bco\label{homcalc}\ \ Let $p$ be an odd prime.\\
(i) $\sigma H_*({\mathcal B}_n(S^2 );\bbf_2)\cong
\bbf_2 [ \iota_{(3,1)},
f_{(5,2)}, f_{(9,4)},\ldots, f_{(2^{i}+1,2^i)},\ldots ]_n$.\\
(ii) $\sigma H_*({\mathcal B}_n(S^2 );\bbf_p)$ is given by
$$\bigoplus_{r+s=n}
E[ \iota_{(3,1)}, h_{(2p+1,p)},\ldots, h_{(2p^i+1,p^i)},\ldots ]_r
\tensor \bbf_p [b_{(2p+2,p)},\ldots ,b_{(2p^i+2,p^i)},\ldots ]_s$$
(iii) For general spheres and $p>n>1$, $\tilde H_*({\mathcal B}_n(S^k);\bbf_p
)$ has a single non-trivial group $\bbf_p$ for $*=n(k+1)-1$ if $k$
is odd, and is trivial if $k$ is even. In particular and for
$p>n>1$, ${\mathcal B}_n(S^{2k})$ has no torsion of order $p$ in its
homology.\eco

\bex The mod-$p$ version of Theorem \ref{jim} shows that
$H^*(\sp{n}S^k,\bbf_p)$ for $p>n$ is a truncated polynomial algebra
$\bbf_p[u]/u^{m+1}$ where $m=1$ or $n$ according to whether $k$ is
odd or even. This fact is recorded in \cite{ucci}, \S1. The element
$u$ pulls back to a generator of $H^{k}S^k$ under the inclusion
$S^k\hookrightarrow\sp{n}S^k$ and so it is of bidegree $(k,1)$. In
particular and for $p>n>1$, $H_*(\bsp{n}S^k;\bbf_p )$ has a single
non-trivial (reduced) homology group for $*=nk$ and $k$ even, and is
entirely trivial if $k$ is odd. This implies corollary \ref{homcalc}
(iii). Note the consistency of this corollary with corollary \ref{topdim}. \eex

\bex As indicated in corollary \ref{homcalc},
$\sigma H_*({\mathcal B}_n(S^2 );\bbf_2)$ consist of elements of filtration
degree $n$ in $\bbf_2 [ \iota_{(3,1)},
f_{(5,2)}, f_{(9,4)},\ldots, f_{(2^{i}+1,2^i)},\ldots ]$.
When $n=2$, there are only two generators:
$f_{(5,2)}$ and $\iota^2_{(3,1)}$ (of bidegree $(6,2)$) so that
$$\tilde H_*({\mathcal B}_2(S^2);\bbf_2) = \begin{cases}\bbf_2, *=4\\
\bbf_2, *=5
\end{cases}$$ and is zero otherwise. For $p$ odd,
$\tilde H_*({\mathcal B}_2(S^2);\bbf_p)$ is trivial (i.e no non trivial
class in filtration degree $2$ in this case since $\iota^2_{(3,1)} = 0$ with odd primes)
and this is consistent with ${\mathcal B}_2(S^2)\simeq\Sigma^{3}\bbr
P^2$  as shown encore in (\ref{jmea}). \eex

%********************************************************************

\addcontentsline{toc}{section}{Bibliography}
\bibliography{biblio}
\bibliographystyle{plain[10pt]}

\enddocument